\documentclass{article}
\pdfoutput=1

\usepackage{arxiv}

\usepackage[utf8]{inputenc} 
\usepackage[T1]{fontenc}    
\usepackage{hyperref}       
\usepackage{url}            
\usepackage{booktabs}       

\usepackage{amsmath,amssymb, amsfonts}       
\usepackage{nicefrac}       
\usepackage{microtype}      
\usepackage{graphicx}
\usepackage{subcaption}
\usepackage[square,numbers]{natbib}
\usepackage{doi}
\usepackage{physics}
\usepackage{bm}
\usepackage{siunitx}
\usepackage{todonotes}
\usepackage[frozencache,cachedir=.]{minted}

\DeclareMathOperator{\fl}{fl}
\DeclareMathOperator{\sign}{sign}
\DeclareMathOperator{\diag}{diag}
\DeclareMathOperator{\vect}{vec}

\newcommand{\trp}{\raisebox{0.2ex}{$\scriptscriptstyle\!\top$}}

\title{Symbolic spectral decomposition of 3x3 matrices}

\author{%
	Michal HABERA \\
	Department of Engineering\\
	University of Luxembourg\\
	\texttt{michal.habera@uni.lu} \\
    \And
    Andreas ZILIAN \\
	Department of Engineering\\
	University of Luxembourg\\
	\texttt{andreas.zilian@uni.lu} \\
}

\renewcommand{\headeright}{Preprint}
\renewcommand{\undertitle}{Preprint}

\hypersetup{
pdftitle={Symbolic spectral decomposition of 3x3 matrices},
pdfauthor={M.~Habera, A.~Zilian},
}

\begin{document}
\maketitle

\begin{abstract}
    Spectral decomposition of matrices is a recurring and important task in applied mathematics, 
    physics and engineering.
    Many application problems require the consideration of matrices of size three with spectral 
    decomposition over the real numbers.
    If the functional dependence of the spectral decomposition on the matrix elements has to be preserved, 
    then closed-form solution approaches must be considered.
	Existing closed-form expressions are based on the use of principal matrix invariants which suffer from 
	a number of deficiencies when evaluated in the framework of finite precision arithmetic.
	This paper introduces an alternative form for the computation of the involved matrix invariants 
	(in particular the discriminant) in terms of sum-of-products expressions as function of the matrix 
	elements.
	We prove and demonstrate by numerical examples that this alternative approach leads to increased
	floating point accuracy, especially in all important limit cases (e.g. eigenvalue multiplicity).
	It is believed that the combination of symbolic algorithms with the accuracy improvements
	presented in this paper can serve as a powerful building block for many engineering tasks.
\end{abstract}


\section{Introduction}

Spectral decomposition of real-valued matrices (eigendecomposition) is a task of utmost 
importance for mathematicians, physicists and engineers. 
Specifically, the decomposition of $3 \times 3$ matrices plays a central role in three 
dimensional space as it is characteristic to many real-world application contexts. 
This special spectral decomposition problem was studied for centuries, especially due to 
its close relation to the roots of the cubic equation. 
Closed-form solutions found increasingly more use with the advent of computers and powerful 
symbolic systems as Mathematica \citep{wolfram1991mathematica} or SymPy \citep{sympy}. 
The symbolic approach represents an efficient tool for computation of eigenvalues and 
eigenvectors while preserving their functional dependence on the matrix elements. 
It has the potential to be exploited together with automatic differentiation (AD) in 
problems where derivatives of spectral decomposition are required 
(e.g. non-linear problems in principal space, sensitivity analysis).

Unfortunately, when results of such eigendecomposition are implemented and evaluated in 
computer software with finite precision arithmetic, not all closed-form approaches are 
equivalent%
\footnote{%
    It might sound confusing to talk about a symbolic algorithm having finite precision issues. 
    What is meant here is that when a symbolic expression gets evaluated for finite precision 
    inputs the result will contain rounding error. 
    It is implicitly assumed that the symbolic engine which evaluates expressions does not 
    perform any aggressive simplifications or optimisations on the expression tree.}%
.
Very few existing papers address accuracy and sensitivity of this decomposition in the 
context of a symbolic algorithm, see e.g. \citep{korelc2014closed, hudobivnik2016closed}.
Papers that study finite precision accuracy of closed-form roots to the cubic equation are more 
common, but less so in the context of eigenvalues and eigenvectors, \citep{kopp2008efficient}. 
Heuristic approaches were developed by the engineering community in order to overcome rounding 
issues, but none effectively improves the accuracy of results, 
see \citep{simo1991quasi, miehe1993computation, jeremic2005significance}.

While well-established iterative methods can provide eigenvalues and eigenvectors up to 
very high precision, their inherent nature 
(presence of loops, non-predetermined number of required iterations, stopping criteria and conditionals) 
renders them not suitable for use within symbolic algorithms. 
Common implementations of iterative schemes can be found in LAPACK \citep{anderson1999lapack}, 
Numerical Recipes in C \citep{numericalrecipes} or GNU Scientific Library \citep{galassi2002gnu}.

The aim of this paper is to study rounding errors in closed-form solution to the 
spectral decomposition. 
Alternative -- but mathematically equivalent -- expressions are sought such that 
rounding effects are diminished. 
An ideal set of symbolic expressions must not require the evaluation of loops or taking 
advantage of any finite precision specific tricks, such that the resulting mathematical 
expression is robust and ready for direct use within a symbolic framework.


\section{Spectral decomposition}
\label{sec:spectral}

For a diagonalisable matrix $\mathbf{A} \in \mathbb R^{n \times n}$ the multiplicative 
spectral decomposition over the real numbers
\begin{align}
    \mathbf{A} = \mathbf{U} \mathbf{\Lambda} \mathbf{V} %
    \label{eq:spectralmult}
\end{align}
is given in terms of the real-valued matrices 
$\mathbf{U} = \left[ \mathbf{u}_1, \dots, \mathbf{u}_n \right]$ and 
$\mathbf{V} = \left[ \mathbf{v}_1^{\trp}, \dots, \mathbf{v}_n^{\trp} \right]^{\trp}$ 
(which, respectively, contain the right eigenvectors $\mathbf{u}_k$ column-wise and the 
left eigenvectors $\mathbf{v}_k$ row-wise) and the diagonal matrix $\mathbf{\Lambda}$, 
holding the $n$ real-valued eigenvalues $\lambda_k$.
The (scaled) eigenvectors in $\mathbf{U}$ and $\mathbf{V}$ fulfil the requirement 
$\mathbf{U}\mathbf{V} = \mathbf{I}$.
The equivalent additive real-valued spectral decomposition 
\begin{align}
    \mathbf{A} = \sum\limits_{k=1}^{n} \lambda_k \mathbf{E}_k
    \label{eq:spectraladd}
\end{align}
can be written in terms of the product of the eigenvalue $\lambda_k$ and its associated 
eigenprojector $\mathbf{E}_k = \mathbf{u}_k \otimes \mathbf{v}_k$.
Independent of eigenvalue multiplicities, $n$ distinct eigenprojectors could always be found such that they have then the following properties
\begin{align}
    \sum\limits_{k=1}^n \mathbf{E}_k &= \mathbf{I}, &
    \mathbf{E}_i \mathbf{E}_j &= \delta_{ij} \mathbf{E}_j, &
    \tr{\mathbf{E}_k} &= 1, &
    \det{\mathbf{E}_k} &= 0
    .
    \label{eq:eigenprojectorsproperties}
\end{align}

\subsection{Eigenvalues and discriminant}
\label{ssec:eigenvalues}

The formulation of the eigenvalue problems 
\begin{subequations} \label{eq:eigenvalueproblem}
\begin{align}
    \mathbf{A}\mathbf{E}_k &= \lambda_k\mathbf{E}_k \label{eq:eigenvalueproblemright} \\
    \mathbf{A}^{\trp}\mathbf{E}_k^{\trp} &= \lambda_k\mathbf{E}_k^{\trp} 
    \label{eq:eigenvalueproblemleft}
\end{align}
\end{subequations}
(or, alternatively, 
$\mathbf{A}\mathbf{U} = \mathbf{U}\mathbf{\Lambda}$ and 
$\mathbf{V}\mathbf{A} = \mathbf{\Lambda}\mathbf{V}$) leads to the characteristic polynomial
\begin{align}
    P_\mathbf{A} (\lambda)
    = \det(\lambda\mathbf{I} - \mathbf{A})
    = \det(\lambda\mathbf{I} - \mathbf{A}^{\trp})
    = \prod\limits_{k=1}^n (\lambda - \lambda_k)
    \label{eq:charpoly}
\end{align}
of matrix $\mathbf{A}$.
The discriminant of the characteristic polynomial $P_\mathbf{A}$ is defined as 
the product of the squared distances of the roots
\begin{align}
    \Delta
    &:=
    \prod\limits_{i < j}^{n} (\lambda_i - \lambda_j)^2
    \label{eq:discriminant}
\end{align}
and is zero in the case of repeated eigenvalues.
The discriminant associated with matrix $\mathbf{A}$ is a function of the matrix elements and 
it has been shown by Parlett \cite{parlett2002discriminant} that the discriminant can be expressed
as the determinant of a symmetric matrix
\begin{align}\label{eq:disc-as-det}
    \Delta = \det{\mathbf{B}} = \det{\mathbf{X}\mathbf{Y}}
\end{align}
with elements 
$B_{ij} = \tr{\mathbf{A}^{i+j-2}} 
        = \mathbf{A}^{i-1} : (\mathbf{A}^{j-1})^{\trp} 
        = \vect^{\trp}\left\{\mathbf{A}^{i-1}\right\} \vect\left\{(\mathbf{A}^{j-1})^{\trp}\right\}$ 
for $1 \leq i, j \leq n$
and a factorisation into the $n \times n^2$ matrix $\mathbf{X}$ and the $n^2 \times n$ matrix $\mathbf{Y}$,
which are -- similar to the Vandermonde matrix -- constructed from powers of $\mathbf{A}$ 
as follows (the $\mathrm{vec}$ operator represents column-stacking)
\begin{align}
    \mathbf{X} &=
    \begin{bmatrix}
    \vect^{\trp}\left\{\mathbf{A}^0\right\} \\
    \vdots \\
    \vect^{\trp}\left\{\mathbf{A}^{n-1}\right\}
    \end{bmatrix},
    &
    \mathbf{Y} &=
    \begin{bmatrix}
    \vect\left\{(\mathbf{A}^0)^{\trp}\right\} &
    \hdots &
    \vect\left\{(\mathbf{A}^{n-1})^{\trp}\right\}
    \end{bmatrix}.
\end{align}
Using the Cauchy-Binet formula for the determinant of a matrix product, the discriminant
is identified as a \textit{sum-of-products}
\begin{align}
    \Delta =\sum\limits_m \det{\mathbf{X}_m} \det{\mathbf{Y}_m} = \mathbf{x}^{\trp} \mathbf{y}
    \label{eq:sop-discriminant}
\end{align}
where the determinants of the $\binom{n^2}{n}$ minors, $\det{\mathbf{X}_m}$ and $\det{\mathbf{Y}_m}$, 
are collected in the vectors $\mathbf{x}$ and $\mathbf{y}$, respectively.
The expansion \eqref{eq:sop-discriminant} has $\binom{n^2}{n} - \binom{n^2 - n}{n}$ nonzero terms
of which several may occur repeatedly.
Thus, a condensed form of the sum-of-products discriminant representation can be established
\begin{align}
    \Delta = \mathbf{\bar x}^{\trp} \mathbf{D} \mathbf{\bar y}
    \label{eq:sop-discriminant-condensed}
\end{align}
in which $\mathbf{\bar x}$ and $\mathbf{\bar y}$ contain only unique factors and 
the diagonal matrix $\mathbf{D}$ holds the respective product multiplier.

\subsubsection{Sub-discriminants}

The expression for the discriminant based on the determinant of matrix $\mathbf B$ suggests a  
generalisation into certain invariants (in this paper called \textit{sub-discriminants}). 
The matrix $\mathbf B = \mathbf X \mathbf Y$ is factored into two rectangular matrices which represent 
column and row-stacked powers $\mathbf A^0, \mathbf A^1, \ldots, \mathbf A^{n-1}$. 
The sub-discriminant $\Delta_{kl}$ corresponding to multi-indices $k = (k_0, k_1, \ldots, k_r)$ 
and $l = (l_0, l_1, \ldots, l_s)$ is defined based on subsets $k$ and $l$ of powers, i.e.
\begin{align}
    \Delta_{kl} := \det(\mathbf B_{kl}) = \det(\mathbf X_k \mathbf Y_l)
    \label{eq:sop-subdiscriminant}
\end{align}
where
\begin{align}
    \mathbf{X}_k &=
    \begin{bmatrix}
    \vect^{\trp}\left\{\mathbf{A}^{k_0}\right\} \\
    \vdots \\
    \vect^{\trp}\left\{\mathbf{A}^{k_r}\right\}
    \end{bmatrix},
    &
    \mathbf{Y}_l &=
    \begin{bmatrix}
    \vect\left\{(\mathbf{A}^{l_0})^{\trp}\right\} &
    \hdots &
    \vect\left\{(\mathbf{A}^{l_s})^{\trp}\right\}
    \end{bmatrix}.
\end{align}
It is important to note, that all sub-discriminants $\Delta_{kl}$ are proper invariants of matrix $\mathbf A$ (they are invariant under similarity transformations). 
Moreover, each component of the matrix $\mathbf B_{kl}$ is an invariant, since 
$(\mathbf B_{kl})_{ij} = \tr{\mathbf A^{i+j-2}}$ for all $i \in k$ and all $j \in l$.

Following its definition simple identities are identified,
\begin{align}
    \Delta_{(0)(1)} &= \Delta_{(1)(0)} = \tr{\mathbf A},\\
    \Delta_{(0, 1, \ldots, n-1)(0, 1, \ldots, n-1)} &= \Delta.
\end{align}

\subsection{Eigenprojectors}
\label{ssec:eigenprojectors}

The $n$ eigenvalues -- as roots $\lambda_k$ of the characteristic polynomial -- are 
nonlinear functions of the $n^2$ elements of $\mathbf{A}$.
A relation for the eigenprojector $\mathbf{E}_k$ is obtained by forming the inner product of
the total differential of \eqref{eq:eigenvalueproblemleft} with $\mathbf{E}_k$ and by using
equation \eqref{eq:eigenvalueproblemright} together with identities 
\eqref{eq:eigenprojectorsproperties} as follows
\begin{subequations}
\begin{align}
    0 &=
    \mathbf{E}_k : \mathrm{d} \left[ (\lambda_k \mathbf{I} - \mathbf{A}^{\trp}) \mathbf{E}_k^{\trp} \right]
    \\
    &=
    \mathbf{E}_k : 
        \left[ 
            (\mathrm{d} \lambda_k \mathbf{I} - \mathrm{d} \mathbf{A}^{\trp}) \mathbf{E}_k^{\trp} 
            +
            (\lambda_k \mathbf{I} - \mathbf{A}^{\trp}) \mathrm{d} \mathbf{E}_k^{\trp} 
        \right]
    \\
    &=
    (\mathrm{d} \lambda_k \mathbf{I} - \mathrm{d} \mathbf{A}) \mathbf{E}_k : \mathbf{E}_k^{\trp} 
    +
    (\lambda_k \mathbf{I} - \mathbf{A})  \mathbf{E}_k : \mathrm{d} \mathbf{E}_k^{\trp} 
    \\
    &=
    \mathrm{d} \lambda_k \tr{\mathbf{E}_k \mathbf{E}_k} - \tr{\mathrm{d} \mathbf{A}\mathbf{E}_k\mathbf{E}_k} 
    +
    \mathbf{0} : \mathrm{d} \mathbf{E}_k^{\trp} 
    \\
    &=
    \mathrm{d} \lambda_k - \mathbf{E}_k^{\trp} : \mathrm{d} \mathbf{A}
    \\
    &= 
    \left( \frac{\partial \lambda_k}{\partial \mathbf{A}} - \mathbf{E}^{\trp} \right) : \mathrm{d} \mathbf{A}
\end{align}
\end{subequations}
from which, for arbitrary $\mathrm{d} \mathbf{A}$, one extracts an expression for the eigenprojector
\begin{align}
    \mathbf{E}_k^{\trp} = \frac{\partial \lambda_k}{\partial \mathbf{A}}
    \label{eq:eigenprojector}
\end{align}
in terms of differentiation through the eigenvalue.
If $\lambda_k = \lambda_k (\mathbf{A})$ is available as differentiable symbolic expression 
(as it is the object of this paper),
then Equation \eqref{eq:eigenprojector} presents advantages over the direct computation of 
the Frobenius covariants as Lagrange interpolants
\begin{align}
    \mathbf{E}_k = \prod\limits_{i=1, i \neq k}^{n} \frac{\mathbf{A} - \lambda_i \mathbf{I}}{\lambda_k - \lambda_i}
\end{align}
obtained as a consequence of Equations \eqref{eq:spectraladd} and \eqref{eq:eigenprojectorsproperties} 
for strictly distinct eigenvalues $\lambda_k$. 
Alternatively, the eigenprojectors can be computed from
\begin{align}\label{eq:eigenproj-from-uv}
    \mathbf E_k = \mathbf U \diag\left\{ \mathbf e_k \right\} \mathbf V
\end{align}
where $\mathbf e_k$ is a canonical unit vector 
(having $1$ in the $k$-th component and $0$ elsewhere) 
and the $\diag$ operator mapping vectors into diagonal matrices. In this case matrices $\mathbf U$ and $\mathbf V$ must be otherwise known, so this approach is useful only for testing purposes.

\section[Spectral decomposition of 3-by-3 matrices]{Spectral decomposition of $3 \times 3$ matrices}
\label{sec:case3by3}

In the following we consider the case of a $3 \times 3$ matrix $\mathbf{A}$ 
with real-valued elements $A_{ij}$ and $1 \leq i, j \leq 3$ under the (above stated) assumption 
that $\mathbf{A}$ is diagonalisable over the real numbers.

\subsection{Eigenvalues}
\label{ssec:eigenvalues3}

The characteristic polynomial of a $3 \times 3$ matrix $\mathbf{A}$
\begin{align}
    P_\mathbf{A} (\lambda)
    = \lambda^3 - I_1(\mathbf{A}) \lambda^2 + I_2(\mathbf{A}) \lambda - I_3(\mathbf{A})
    = (\lambda - \lambda_1)(\lambda - \lambda_2)(\lambda - \lambda_3)
    \label{eq:charpoly3}
\end{align}
can be expressed in terms of the principal invariants of $\mathbf{A}$
\begin{align}
    I_1(\mathbf{A}) &= \tr{\mathbf{A}}, &
    I_2(\mathbf{A}) &= \frac12 \left( (\tr{\mathbf{A}})^2 - \tr{\mathbf{A}^2} \right), &
    I_3(\mathbf{A}) &= \det{\mathbf{A}}
    \label{eq:invariants}
\end{align}
or in terms of its three roots $\lambda_1$, $\lambda_2$ and $\lambda_3$.

\subsubsection{Roots of the characteristic polynomial}
\label{sssec:charpolyroots}

In order to facilitate the identification of the roots $\lambda_k$ of the cubic equation 
$P_\mathbf{A} (\lambda) = 0$
the substitution
\begin{align}
    \lambda = \mu + r \, \cos \varphi
    \label{eq:charpolysubst}
\end{align}
is introduced, with $\mu, r,\varphi \in \mathbb{R}$ to be determined.
This leads, in a first step, to
\begin{align}
    P_\mathbf{A} (\mu, r, \varphi) 
    = \left( r^3 \right) \cos^3 \varphi
    + r^2 \left( 3\mu - I_1 \right) \cos^2 \varphi
    + r \left( 3\mu^2 - 2 I_1 \mu + I_2 \right) \cos \varphi
    + \left( \mu^3 - I_1 \mu^2 + I_2 \mu - I_3 \right)
    \nonumber
\end{align}
and is transformed to
\begin{align}
    P_\mathbf{A} (\mu, r, \varphi) 
    = \frac{r^3}{4} \cos 3\varphi
    + r^2 \left( 3\mu - I_1 \right) \cos^2 \varphi
    + \frac{r}{4} \left( 12\mu^2 - 8 I_1 \mu + 4 I_2 + 3 r^2 \right) \cos \varphi
    + \left( \mu^3 - I_1 \mu^2 + I_2 \mu - I_3 \right)
    \nonumber
\end{align}
using the triple-angle trigonometric identity $4 \cos^3 \varphi = \cos 3\varphi - 3 \cos \varphi$.
In the above expression for $P_\mathbf{A}$ the factors to $\cos^2 \varphi$ and $\cos \varphi$ become
zero for suitable choices of $\mu$ and $r$.
The non-trivial solution of the associated system of equations
\begin{align}
    r^2 \left( 3\mu - I_1 \right) &= 0 & 
    \frac{r}{4} \left( 12\mu^2 - 8 I_1 \mu + 4 I_2 + 3 r^2 \right) &= 0 
\end{align}
is given by
\begin{align}
    \mu &= \frac13 I_1
    & 
    r &= \pm \frac23 \sqrt{I_1^2 - 3 I_2}
\end{align}
and allows to eliminate $\mu$ and $r$ in $P_\mathbf{A}(\mu, r, \varphi)$. This provides
\begin{align}
    P_\mathbf{A} (\varphi) 
    &= \pm \frac{2}{27} \left( I_1^2 - 3 I_2 \right)^\frac32 \cos 3\varphi 
     + \left( -\frac{2}{27} I_1^3 + \frac13 I_1 I_2 - I_3 \right) 
    \nonumber
\end{align}
and with
\begin{align}
    0
    \overset{!}{=}
    \pm 2 \left( I_1^2 - 3 I_2 \right)^\frac32 \cos 3\varphi + 
       \left( -2 I_1^3 + 9 I_1 I_2 - 27 I_3 \right)
\end{align}
the relation to determine the cosine of the triple angle
\begin{align}
    \cos 3\varphi
    &=
    \pm \frac12 \left( 2 I_1^3 - 9 I_1 I_2 + 27 I_3 \right) \left( I_1^2 - 3 I_2 \right)^{-\frac32}
    \label{eq:triple-angle}
\end{align}
with $-1 \leq \cos 3\varphi \leq +1$.

The three eigenvalues of matrix $\mathbf{A}$
\begin{align}\label{eq:eigenvalues3}
    \lambda_k
    &=
    \frac13 \left[ I_1 + 2 \sqrt{I_1^2 - 3 I_2} \cos (\varphi + \frac{2 \pi}{3} k) \right],
    \qquad
    k = 1, 2, 3
\end{align}
are then highly nonlinear functions 
of the principal matrix invariants.

\subsubsection{Discriminant and sub-discriminants}
\label{sssec:charpolydisc}

The discriminant of the (cubic) polynomial can be expressed in terms of its coefficients (referred to as ``naive'' expression in this paper)
\begin{subequations}
\begin{align}
    \Delta
    &=
    18 I_1 I_2 I_3 + I_1^2 I_2^2 - 4 I_1^3 I_3 - 4 I_2^3 - 27 I_3^2 
    \label{eq:naive-discriminant}
    \\
    &=
    \frac{1}{27}
        \left[
            4 \left( I_1^2 - 3 I_2 \right)^3 - \left( 2 I_1^3 - 9 I_1 I_2 + 27 I_3 \right)^2
        \right]
    =
    \frac{1}{27}
        \left(
            4 \Delta_p^3 - \Delta_q^2
        \right)
    \label{eq:dpdq-discriminant}
\end{align}
\end{subequations}
where invariants $\Delta_p$ and $\Delta_q$ are defined as
\begin{subequations} \label{eq:def-dpdq}
\begin{align} 
    \Delta_p & := I_1^2 - 3 I_2, \label{eq:naive-deltap} \\
    \Delta_q & := 2 I_1^3 - 9 I_1 I_2 + 27 I_3. \label{eq:naive-deltaq}
\end{align}
\end{subequations} 

Alternatively, for $n = 3$ the condensed sum-of-products representation of the discriminant, 
see Equation \eqref{eq:sop-discriminant-condensed}, is given by 14 products with
\begin{subequations} \label{eq:sop-discriminant-condensed3}
\begin{align}
    \mathbf{\bar x} &= 
    \scalebox{0.71}{$
    \begin{bmatrix}
    A_{12} A_{23} A_{31} - A_{13} A_{21} A_{32} \\
    A_{12}^{2} A_{23} - A_{12} A_{13} A_{22} + A_{12} A_{13} A_{33} - A_{13}^{2} A_{32} \\
    A_{11} A_{12} A_{32} - A_{12}^{2} A_{31} - A_{12} A_{32} A_{33} + A_{13} A_{32}^{2} \\
    A_{11} A_{13} A_{23} + A_{12} A_{23}^{2} - A_{13}^{2} A_{21} - A_{13} A_{22} A_{23} \\
    A_{11} A_{12} A_{23} - A_{12} A_{13} A_{21} - A_{12} A_{23} A_{33} + A_{13} A_{23} A_{32} \\
    A_{11} A_{13} A_{32} - A_{12} A_{13} A_{31} + A_{12} A_{23} A_{32} - A_{13} A_{22} A_{32} \\
    A_{12} A_{21} A_{23} - A_{13} A_{21} A_{22} + A_{13} A_{21} A_{33} - A_{13} A_{23} A_{31} \\
    A_{11}^{2} A_{23} - A_{11} A_{13} A_{21} - A_{11} A_{22} A_{23} - A_{11} A_{23} A_{33} + A_{12} A_{21} A_{23} + A_{13} A_{21} A_{33} + A_{22} A_{23} A_{33} - A_{23}^{2} A_{32} \\
    A_{11}^{2} A_{23} - A_{11} A_{13} A_{21} - A_{11} A_{22} A_{23} - A_{11} A_{23} A_{33} + A_{13} A_{21} A_{22} + A_{13} A_{23} A_{31} + A_{22} A_{23} A_{33} - A_{23}^{2} A_{32} \\
    A_{11} A_{12} A_{22} - A_{11} A_{12} A_{33} - A_{12}^{2} A_{21} + A_{12} A_{13} A_{31} - A_{12} A_{22} A_{33} + A_{12} A_{33}^{2} + A_{13} A_{22} A_{32} - A_{13} A_{32} A_{33} \\
    A_{11} A_{12} A_{22} - A_{11} A_{12} A_{33} + A_{11} A_{13} A_{32} - A_{12}^{2} A_{21} - A_{12} A_{22} A_{33} + A_{12} A_{23} A_{32} + A_{12} A_{33}^{2} - A_{13} A_{32} A_{33} \\
    A_{11} A_{12} A_{23} - A_{11} A_{13} A_{22} + A_{11} A_{13} A_{33} - A_{12} A_{22} A_{23} - A_{13}^{2} A_{31} + A_{13} A_{22}^{2} - A_{13} A_{22} A_{33} + A_{13} A_{23} A_{32} \\
    A_{11} A_{13} A_{22} - A_{11} A_{13} A_{33} - A_{12} A_{13} A_{21} + A_{12} A_{22} A_{23} - A_{12} A_{23} A_{33} + A_{13}^{2} A_{31} - A_{13} A_{22}^{2} + A_{13} A_{22} A_{33} \\
    A_{11}^{2} (A_{22} - A_{33}) + A_{22}^{2} (A_{33} - A_{11}) + A_{33}^{2} (A_{11} - A_{22})
    + A_{11} (A_{13} A_{31}- A_{12} A_{21}) + A_{22} (A_{12} A_{21} - A_{23} A_{32}) + A_{33} (A_{23} A_{32} - A_{13} A_{31})
    \end{bmatrix}
    $}, \\
    \mathbf{\bar y} &= 
    \scalebox{0.71}{$
    \begin{bmatrix}
    A_{13} A_{21} A_{32} - A_{12} A_{23} A_{31} \\
    A_{21}^{2} A_{32} - A_{21} A_{22} A_{31} + A_{21} A_{31} A_{33} - A_{23} A_{31}^{2} \\
    A_{11} A_{21} A_{23} - A_{13} A_{21}^{2} - A_{21} A_{23} A_{33} + A_{23}^{2} A_{31} \\
    A_{11} A_{31} A_{32} - A_{12} A_{31}^{2} + A_{21} A_{32}^{2} - A_{22} A_{31} A_{32} \\
    A_{11} A_{21} A_{32} - A_{12} A_{21} A_{31} - A_{21} A_{32} A_{33} + A_{23} A_{31} A_{32} \\
    A_{11} A_{23} A_{31} - A_{13} A_{21} A_{31} + A_{21} A_{23} A_{32} - A_{22} A_{23} A_{31} \\
    A_{12} A_{21} A_{32} - A_{12} A_{22} A_{31} + A_{12} A_{31} A_{33} - A_{13} A_{31} A_{32} \\
    A_{11}^{2} A_{32} - A_{11} A_{12} A_{31} - A_{11} A_{22} A_{32} - A_{11} A_{32} A_{33} + A_{12} A_{21} A_{32} + A_{12} A_{31} A_{33} + A_{22} A_{32} A_{33} - A_{23} A_{32}^{2} \\
    A_{11}^{2} A_{32} - A_{11} A_{12} A_{31} - A_{11} A_{22} A_{32} - A_{11} A_{32} A_{33} + A_{12} A_{22} A_{31} + A_{13} A_{31} A_{32} + A_{22} A_{32} A_{33} - A_{23} A_{32}^{2} \\
    A_{11} A_{21} A_{22} - A_{11} A_{21} A_{33} - A_{12} A_{21}^{2} + A_{13} A_{21} A_{31} - A_{21} A_{22} A_{33} + A_{21} A_{33}^{2} + A_{22} A_{23} A_{31} - A_{23} A_{31} A_{33} \\
    A_{11} A_{21} A_{22} - A_{11} A_{21} A_{33} + A_{11} A_{23} A_{31} - A_{12} A_{21}^{2} - A_{21} A_{22} A_{33} + A_{21} A_{23} A_{32} + A_{21} A_{33}^{2} - A_{23} A_{31} A_{33} \\
    A_{11} A_{21} A_{32} - A_{11} A_{22} A_{31} + A_{11} A_{31} A_{33} - A_{13} A_{31}^{2} - A_{21} A_{22} A_{32} + A_{22}^{2} A_{31} - A_{22} A_{31} A_{33} + A_{23} A_{31} A_{32} \\
    A_{11} A_{22} A_{31} - A_{11} A_{31} A_{33} - A_{12} A_{21} A_{31} + A_{13} A_{31}^{2} + A_{21} A_{22} A_{32} - A_{21} A_{32} A_{33} - A_{22}^{2} A_{31} + A_{22} A_{31} A_{33} \\
    A_{11}^{2} (A_{22} - A_{33}) + A_{22}^{2} (A_{33} - A_{11}) + A_{33}^{2} (A_{11} - A_{22})
    + A_{11} (A_{13} A_{31}- A_{12} A_{21}) + A_{22} (A_{12} A_{21} - A_{23} A_{32}) + A_{33} (A_{23} A_{32} - A_{13} A_{31})
    \end{bmatrix}
    $}, \\
    \mathbf{D} &= \diag\{(9, 6, 6, 6, 8, 8, 8, 2, 2, 2, 2, 2, 2, 1)\}
    .
\end{align}
\end{subequations}
It should be noted that \emph{symmetric} matrices allow for further reduction of the number 
of products and presentation of the discriminant as a \textit{sum-of-squares}, as demonstrated 
by Kummer \cite{kummer1843bemerkungen} (7 squares) and Watson \cite{watson1956some} (5 squares).

With the notion of sub-discriminants it can be shown that
\begin{align}
    \Delta_{(0, 1)(0, 1)} = 2 \Delta_p = 2 (I_1^2 - 3I_2).
\end{align}

It follows from
\begin{align}
    \det{\mathbf B_{(0, 1)(0, 1)}} = \det{\begin{bmatrix}
        3 & \tr{A} \\
        \tr{A} & \tr{A^2}
    \end{bmatrix}} = 3 \tr{A^2} - \tr{A}^2 = 2 \Delta_p.
\end{align}

With similar factorisation into sum-of-products an alternative expression for the $\Delta_p$ 
invariant reads
\begin{subequations} \label{eq:sop-deltap}
\begin{align}
\Delta_p &= \frac12 \mathbf{\bar x}_p^{\trp} \mathbf D_p \mathbf{\bar y}_p, \\
\mathbf{\bar x}_p &= \begin{bmatrix}
    A_{10} \\
    A_{20} \\
    A_{21} \\
    -A_{00} + A_{11} \\
    -A_{00} + A_{22} \\
    -A_{11} + A_{22}
\end{bmatrix}, \quad
\mathbf{\bar y}_p = \begin{bmatrix}
    A_{01} \\
    A_{02} \\
    A_{12} \\
    -A_{00} + A_{11} \\
    -A_{00} + A_{22} \\
    -A_{11} + A_{22}
\end{bmatrix},\\
\mathbf D_p &= \diag\{(6, 6, 6, 1, 1, 1)\}.
\end{align}
\end{subequations}

A similar expression cannot be identified for $\Delta_q$, however, this invariant can
partially be expressed on the basis of sum-of-products as follows
\begin{align} \label{eq:sop-deltaq}
    \Delta_q &= 3 \Delta_{(0,1)(0,2)} - 4 \tr{\mathbf A} \Delta_p.
\end{align}

\section{Finite precision and rounding errors}
\label{sec:finite}

The above formulas for the computation of eigenvalues and eigenprojectors are not all equivalent 
when finite precision floating point arithmetic is involved. 
This is, of course, the case when these expressions are implemented in a computer program. 
The discrepancies will be explained and addressed below.

Machine epsilon, $\varepsilon$, is here defined as the difference between 1 and the next larger 
floating point number. 
Let $b$ be the radix (base) and $p$ a precision of a floating point number representation, then 
$\varepsilon = b^{-(p-1)}$. 
For example, the standard double precision floating point representation used in this paper has 
$\varepsilon = \SI{2.22e-16}{}$. 
Only the magnitude of the rounding error is considered in this paper, so there is no distinction 
made between units in the last place and relative error in terms of machine epsilon.

\subsection{Rounding errors in principal invariants}

An inherent property of the analytical approach to spectral decomposition follows as a 
consequence of its use of matrix invariants $I_1, I_2$ and $I_3$.
As noted in \citep{kopp2008efficient}, this approach is not suitable when the distances between 
eigenvalues are large.  As a simple example, consider a matrix with 
$\lambda_1 = 1, \lambda_2 = 2$ and $\lambda_3 = \SI{e+17}{}$. 
Any computation which then makes use of matrix trace $I_1 = \lambda_1 + \lambda_2 + \lambda_3$ 
will suffer from rounding error (in standard double precision), to the effect that the presence 
of the two small eigenvalues is not measurable in the matrix trace.

Due to non-associativity of the floating point arithmetic the case with eigenvalues 
$\lambda_1 = -\SI{e+17}{}, \lambda_2 = 1$ and $\lambda_3 = \SI{e+17}{}$ would also introduce 
large relative error into the computation of the trace $I_1 = \fl(\fl(\lambda_1 + \lambda_2) + \lambda_3)$. 
Here, $\fl(x)$ denotes floating point representation of number $x$. 
Such issues could be solved using techniques as the Kahan summation algorithm 
(or the improved Kahan-Babu\v{s}ka algorithm \citep{neumaier1973rundungsfehleranalyse}). 
Unfortunately, tricks based on benign cancellation are not applicable when developing a symbolic 
algorithm, since the actual execution order of floating point operations is not strictly preserved.

In \citep{kopp2008efficient} the alternative use of iterative methods is recommended (Jacobi/QR and QL). These methods, however, are not relevant in the context of symbolic computation.

\subsection{Catastrophic cancellation in the discriminant}

Catastrophic cancellation is a numerical phenomenon where the subtraction of the good approximation 
to two close numbers results in a bad approximation in the result. 
The formula for discriminant based on principal invariants (Equation \eqref{eq:naive-discriminant} 
and \eqref{eq:dpdq-discriminant}) suffers from this effect.

Let $\mathbf A \in \mathbb{R}^{3 \times 3}$ be a diagonalisable traceless matrix with real 
eigenvalues $\lambda_1 = 1 + \delta, \lambda_2 = 1$ and $\lambda_3 = -2 - \delta$. 
The parameter $\delta$ is a distance between two largest eigenvalues and for 
$\delta \longrightarrow 0$ there is $\Delta \longrightarrow 0$.
For this matrix
\begin{align}
    I_1 = 0, \qquad I_2 = -3 - 3\delta - \delta^2, \qquad I_3 = -2 - 3 \delta - \delta^2.
    \nonumber
\end{align}
If $\delta^2$ is smaller than machine epsilon, then the computation and storage of $I_2$ and $I_3$ 
introduces a rounding error, such that
\begin{align}
    \fl(I_2) = -3 - 3 \delta, \qquad \fl(I_3) = -2 - 3 \delta.
    \nonumber
\end{align}
The discriminant computed using formula based on principal invariants \eqref{eq:naive-discriminant} then
leads to
\begin{align}
    \fl(\Delta) &= \fl(-4\fl(I_2)^3) - \fl(27 \fl(I_3)^2) 
    \nonumber\\
    &= \fl(108 + 324 \delta + 324 \delta^2 + 108 \delta^3) - \fl(108 + 324 \delta + 243 \delta^2)
    \nonumber\\
    &= 0. 
    \nonumber
\end{align}
This situation occurs for standard double precision when $\delta < \SI{e-8}{}$. 
In other words, the discriminant computed from the naive expression \eqref{eq:naive-discriminant} 
becomes insignificant when the distance between two eigenvalues is smaller than half precision.
The loss of precision could be avoided if a different expression for the discriminant is used, 
which motivates the following discussion.

Fortunately, the sum-of-products expression for the discriminant 
(Equation \eqref{eq:sop-discriminant-condensed}) does not suffer from catastrophic cancellation 
when the discriminant goes to zero. 
This observation is the foundation for computation of eigenvalues 
(and determination of eigenvalues multiplicity) 
with improved accuracy and is proved in the following.

For the special case of \emph{symmetric} matrix, sum-of-products reduces to sum-of-squares, 
i.e. the vectors $\mathbf{\bar x}$ and $\mathbf{\bar y}$ are equal. 
Assuming that the distance between the two closest eigenvalues is some $\delta << 1$, the 
discriminant is proportional to $\delta^2$ (since the discriminant is the square of eigenvalue 
distances) and every (positive) summand in the scalar product 
$\mathbf{\bar x}^{\trp} \mathbf D \mathbf{\bar x}$ is thus either zero or a small number 
proportional to $\delta^2$. 
Every non-zero element in the vector $\mathbf{\bar x}$ is therefore proportional to $\delta$ and 
could be computed up to machine precision. Finally, the sum-of-squares of small numbers of the same 
magnitude does not introduce a significant rounding error.

In the general \emph{non-symmetric} case, it must be shown that each non-zero element in 
$\mathbf{\bar x}$ and $\mathbf{\bar y}$ approaches zero dominated by terms linear in $\delta$. 
First, the limit case for $\delta=0$ (i.e. $\Delta=0$) is proven. 
Equivalently, it is required that
\begin{align}
    \Delta = 0 \Longrightarrow \det{\mathbf{X}_m} = \det{\mathbf{Y}_m} = 0, \quad \forall m,
\end{align}
where $m$ runs through all $n \times n$ square sub-matrices of $\mathbf X$ and $\mathbf Y$. 
This is a stronger statement than $\Delta = 0 \Longrightarrow \det(\mathbf X \mathbf Y) = 0$ 
which follows easily from Equation \eqref{eq:disc-as-det}.

It could be shown, that for a diagonalisable matrix $\mathbf A \in \mathbb R^{n \times n}$, 
which has zero discriminant, the following matrices
\begin{align}
    \mathbf I, \mathbf A, \mathbf A^2, \ldots, \mathbf A^{n-1}
\end{align}
are linearly dependent%
\footnote{%
    Matrix $\mathbf A$ is a root of its minimal polynomial, which is of order $k \leq n$. 
    For a diagonalisable matrix every eigenvalue has its algebraic multiplicity equal to 
    geometric multiplicity. 
    Therefore, a zero discriminant implies that the degree of minimal polynomial is strictly smaller 
    than the matrix size $n$.}%
.
In other words, there exist coefficients $\{c_i\} \in \mathbb R^n$ such that
\begin{align}\label{eq:lin-dependence-powers}
    c_0 \mathbf I + c_1 \mathbf A + \ldots + c_{n-1} \mathbf A^{n-1} = \mathbf 0.
\end{align}
Alternatively, \eqref{eq:lin-dependence-powers} could be written using the $\vect$ operator as
\begin{align}
    c_0 \vect^{\trp} \mathbf I + c_1 \vect^{\trp} \mathbf A + \ldots + c_{n-1} \vect^{\trp} \mathbf A^{n-1} = \mathbf 0.
\end{align}
This equation states that the rows in matrix $\mathbf X$ are linearly dependent, thus the rows 
in all sub-matrices $\mathbf X_m$ are linearly dependent too. 
Hence, $\det\{ \mathbf X_m\} = 0$ for each $m$. 
The same observation applies to the transpose of \eqref{eq:lin-dependence-powers} and columns of 
$\mathbf Y_m$, so that equally, $\det\{\mathbf Y_m\} = 0$.
Moreover, since the discriminant is a continuous function of distances between eigenvalues -- 
follows from its definition \eqref{eq:discriminant} --, if $\delta \longrightarrow 0$ then each 
$\det\{\mathbf X_m\} \longrightarrow 0$ and $\det\{\mathbf Y_m\} \longrightarrow 0$. Finally, since
discriminant is proportional to $\delta^2$, each non-zero element in $\mathbf{\bar x}$
and $\mathbf{\bar y}$ is proportional to $\delta$.

\paragraph{Note:} %
Catastrophic cancellation occurs also in \eqref{eq:eigenvalues3} for 
$(\varphi \rightarrow 0, I_1 > 0)$ or $(\varphi \rightarrow \pi, I_1 < 0)$. 
If the distance between eigenvalues $\delta$ becomes large, Equation \eqref{eq:eigenvalues3} 
results in the subtraction of two large numbers in order to produce the small eigenvalue 
$\lambda_2 = 1$ (for the example introduced at the beginning of this subsection).

\subsection[Catastrophic cancellation in the delta-p invariant]{Catastrophic cancellation in the $\Delta_p$ invariant}

The invariant $\Delta_p$ plays a foremost role in the computation of eigenvalues. 
For $\mathbf A \in \mathbb{R}^{3 \times 3}$ being a diagonalisable matrix with eigenvalues 
$\lambda_1 = 1, \lambda_2 = 1$ and $\lambda_3 = 1 + \delta$ one obtains $\Delta_p = \delta^2$. 
For this matrix
\begin{align}
    I_1 = 3 - \delta, \quad I_2 = 3 - 2 \delta.
    \nonumber
\end{align}
Similarly to arguments for cancellation in the discriminant, when $I_1^2$ is computed and 
stored for $\delta^2$ smaller than machine epsilon, one observes that
\begin{align}
    \fl(I_1^2) &= \fl(9 - 6\delta + \delta^2) = 9 - 6\delta,
    \nonumber\\
    \fl(\Delta_p) &= \fl(9 - 6\delta - 3\fl(3 - 2\delta)) = 0.
    \nonumber
\end{align}
Thankfully, the sum-of-products expression for $\Delta_p$, see Equation \eqref{eq:sop-deltap}, 
does not suffer from catastrophic cancellation by similar arguments as used in the
proof for the discriminant. It could be shown that for $\Delta_p$ proportional to $\delta^2$
every non-zero term in vectors $\mathbf{\bar x}_p$ and $\mathbf{\bar y}_p$ is proportional to
$\delta$. It is a consequence of the linear dependence of matrices $\mathbf I$ and $\mathbf A$,
which follows from that fact that degree of minimal polynomial $k = 1$ in the case of $3 \times 3$
matrix with $\Delta_p = 0$.

\subsection{Rounding error in the triple angle}

The evaluation of the triple angle $\varphi$ in Equation \eqref{eq:triple-angle} suffers from 
dangerous rounding errors in the vicinity of $\Delta = 0$. 
In order to demonstrate these effects, let again consider a traceless matrix with eigenvalues 
$\lambda_1 = 1 + \delta, \lambda_2 = 1$ and $\lambda_3 = -2 - \delta$. 
Direct inversion of \eqref{eq:triple-angle} gives
\begin{align}\label{eq:triple-angle-arccos}
    \varphi = \frac13 \arccos \left( \frac12 \Delta_q \Delta_p^{-3/2} \right)
\end{align}
and series expansion of the inverse cosine argument around the critical point $\Delta = 0$ 
(i.e. $\delta = 0$) for the chosen example matrix leads to
\begin{align}
    \frac12 \Delta_q \Delta_p^{-3/2} = -1 + \frac{3\delta^2}{8} + \mathcal{O}(\delta^3).
    \nonumber
\end{align}
Similar to the discussion above, when the argument is evaluated and stored in floating point representation 
rounding error becomes significant if the eigenvalue distance is smaller than square root of machine epsilon. 
This effect could be avoided in two ways.

The first approach is based on the generalised (Puiseux) series approximation to the inverse cosine
\begin{align}
    \acos \left( 1-x \right) = \sqrt{2x} + \frac{x^{3/2}}{6 \sqrt{2}} + \frac{3 x^{5/2}}{80 \sqrt{2}} + \frac{x^{7/2}}{448 \sqrt{2}} + \mathcal{O}(x^{9/2})
\end{align}
around point $x = 0^+$. 
The presence of the square root in the expansion suggests a possible series approximation of the 
angle which would contain linear terms in the eigenvalue distance $\delta$ 
(or equivalently, terms proportional to $\sqrt{\Delta})$). 
Indeed, the expansion of the angle $\varphi$ around $\Delta = 0$ reads
\begin{align}\label{eq:angle-expansion}
    \varphi = \frac13 \acos \left( \frac{\Delta_q}{ \sqrt{\Delta_q^2 + 27 \Delta}} \right) = \frac13 \left[ \acos \left( \sign \Delta_q \right) + \frac{3 \sqrt{3} \sqrt{\Delta}}{\Delta_q} - \frac{27 \sqrt{3} (\sqrt{\Delta})^3}{\Delta_q^3} + \frac{2187\sqrt{3} (\sqrt{\Delta})^5}{5 \Delta_q^5} + \mathcal{O}(\Delta^3)\right].
\end{align}

In the alternative second approach, a trigonometric identity for $\arccos$ is used to transform 
its argument, such that the evaluation of problematic arguments 
(like $-1 + \delta^2$ in the above example) is avoided. 
If the argument is transformed to the vicinity of zero, then the evaluation does not suffer from 
finite-precision round-off.
Applying the Pythagorean theorem to the right-angled triangle with unit hypotenuse leads to
\begin{align}
    \tan \left( \arccos x \right) = \frac{\sqrt{1 - x^2}}{x}.
\end{align}
The inverse tangent could be applied to this identity. 
However, $\arccos$ maps to $[0, \pi]$, so for $x < 0$ the standard $\arctan$ must be 
shifted by $\pi$. 
Hence,
\begin{align}
    \arccos x = \arctan \left( \frac{\sqrt{1 - x^2}}{x} \right) - \frac{\sign(x) - 1}{2} \pi.
\end{align}
With $x = \Delta_q / \sqrt{\Delta_q^2 + 27 \Delta}$ the angle $\varphi$ is equivalently expressed as
\begin{align}
    \varphi = \frac13 \left[ \arctan \left( \frac{3\sqrt{3} \sqrt{\Delta}}{\Delta_q} \right) - \frac{\sign (\Delta_q) - 1}{2} \pi \right].
    \label{eq:angleatan}
\end{align}
The use of the $\arctan$ identity was already noted in \citep{smith1961eigenvalues} for better 
accuracy, but without thorough explanation. 
In addition, the sum-of-products expression for the discriminant must be employed along the 
$\arctan$ identity in order to obtain a consistently increased accuracy in the eigenvalues.

\subsection{Summary of the rounding errors}

The improved expression for eigenvalues is given in terms of three invariants $I_1, \Delta_p$ 
and $\Delta_q$
\begin{align}\label{eq:eigenvalues-dpdq}
    \lambda_k &= \frac13 \left[ I_1 + 2 \sqrt{\Delta_p} \cos (\varphi + \frac{2 \pi}{3} k) \right], \qquad k = 1, 2, 3,
\end{align}
where Equation \eqref{eq:angleatan} is used to compute $\varphi$ based on $\Delta$ and $\Delta_q$. 
With the use of the sum-of-products expression rounding errors are significantly reduced, 
see \autoref{tab:summary-errors}.

\begin{table}
\caption{Summary of rounding errors in invariants}
\label{tab:summary-errors}
\centering
    \begin{tabular}{lllll}\toprule
         invariant & expression & critical case & eigenvalues & absolute error \\ \midrule
         $I_1$ & \eqref{eq:invariants} & large eigenvalue distances & $\lambda_1 = 1, \lambda_2 = 1, \lambda_3 = \SI{e+17}{}$ & 1 \\ \midrule
         $\Delta$  & "naive", \eqref{eq:naive-discriminant} & $\Delta \longrightarrow 0$ & $\lambda_1 = -1, \lambda_2 = 1, \lambda_3 = 1 + \varepsilon$ & $\varepsilon$\\
          & sum-of-products, \eqref{eq:sop-discriminant-condensed3} & & & $\varepsilon^2$\\ \midrule
         $\Delta_p$ & "naive", \eqref{eq:naive-deltap} & $\Delta_p \longrightarrow 0$ & $\lambda_1 = 1, \lambda_2 = 1, \lambda_3 = 1 + \varepsilon$ & $\varepsilon$\\
          & sum-of-products, \eqref{eq:sop-deltap} & & & $\varepsilon^2$\\ \midrule
         $\Delta_q$ & "naive", \eqref{eq:naive-deltaq} & $\Delta_q \longrightarrow 0$ & $\lambda_1 = 0, \lambda_2 = 1, \lambda_3 = 2$ & $\varepsilon$\\
         \bottomrule
    \end{tabular}
\end{table}

\paragraph{Remark:} %
The authors are not aware of an alternative expression for the invariant $\Delta_q$ which would 
exhibit an error smaller than machine precision in the critical case. 
One could show that $\Delta_q = 3 \Delta_{(0, 1)(0, 2)} - 4 \tr{\mathbf A} \Delta_p$, but this 
formula shows cancellation and larger error too. 
An expression based on the sum of small products would be required. 
There are approaches derived from Kahan's algorithm for the determinant of $2 \times 2$ 
matrices \cite{jeannerod2013further}. 
Unfortunately, these techniques benefit from benign cancellation and fused multiply-add 
instructions, and are thus not suitable in the symbolic scope of this work.

\paragraph{Remark:} %
Errors in \autoref{tab:summary-errors} are included for eigenvalues of order 1. 
A simple scaling argument could be used for rough estimates of the error for a matrix with 
elements of different order. 
The discriminant $\Delta$ is a sixth order polynomial in the matrix elements, thus its error in 
\autoref{tab:summary-errors} could be scaled with $\abs{\max(A_{ij})}^6$. 
Similarly, the $\Delta_p$ invariant as a second order polynomial scales with $\abs{\max(A_{ij})}^2$ 
and $\Delta_q$ with $\abs{\max(A_{ij})}^3$.

\section{Numerical benchmarks}
\label{sec:benchmarks}

In order to test the improved accuracy of eigenvalues (and invariants) computation a non-symmetric 
diagonalisable matrix 
$\mathbf B = \mathbf U \mathbf {\Lambda} \mathbf U^{-1} \in \mathbb R^{3 \times 3}$ is considered.

There are three critical cases where invariants $\Delta, \Delta_p$ or $\Delta_q$ vanish. 
In the following benchmarks a distance parameter $\delta$ is used to approach these limit states. 
It can be modelled on the basis of the diagonal matrix $\mathbf{\Lambda}$, 
see \autoref{tab:critical-cases-bench}.

The similarity transformation matrix $\mathbf U$ does not change the value of eigenvalues 
(or invariants) in infinite precision. 
However, it effectively propagates the distance parameter $\delta$ into the elements of the final 
matrix $\mathbf A$ as well as into the principal invariants. 
For this reason two different non-singular transformation matrices are tested, 
see \autoref{tab:similarity-trans-bench}.

\begin{table}
\caption{Critical cases studied in the numerical benchmarks.}
\label{tab:critical-cases-bench}
\centering
    \begin{tabular}{llll}\toprule
          & $\Delta \longrightarrow 0$ & $\Delta_p \longrightarrow 0$ & $\Delta_q \longrightarrow 0$ \\ \midrule
         $\mathbf{\Lambda}$ & $\diag\{(-1, 1, 1 + \delta)\}$ & $\diag\{(1, 1, 1 + \delta)\}$ & $\diag\{(0, 1, 2 + \delta)\}$ \\
         \bottomrule
    \end{tabular}
\end{table}

\begin{table}
\caption{Similarity transformation matrices $\mathbf U$ used in the numerical benchmarks.}
\label{tab:similarity-trans-bench}
\centering
    \begin{tabular}{ccc}\toprule
          & Case I & Case II \\ \midrule
         $\mathbf{U}$ & $\begin{bmatrix}
             1 & -1 & 1 \\
             1 & 1 & 1 \\
             -1 & -1 & 1
         \end{bmatrix}$ & $\begin{bmatrix}
             1 & 1 & 1 \\
             1 & 0 & 1 \\
             2 & 1 & 2  + \gamma
         \end{bmatrix}$ \\
         \bottomrule
    \end{tabular}
\end{table}

\subsection{Invariants}

\begin{figure}
    \centering
    \begin{subfigure}[b]{0.45\textwidth}
    \centering
    \includegraphics[width=\textwidth]{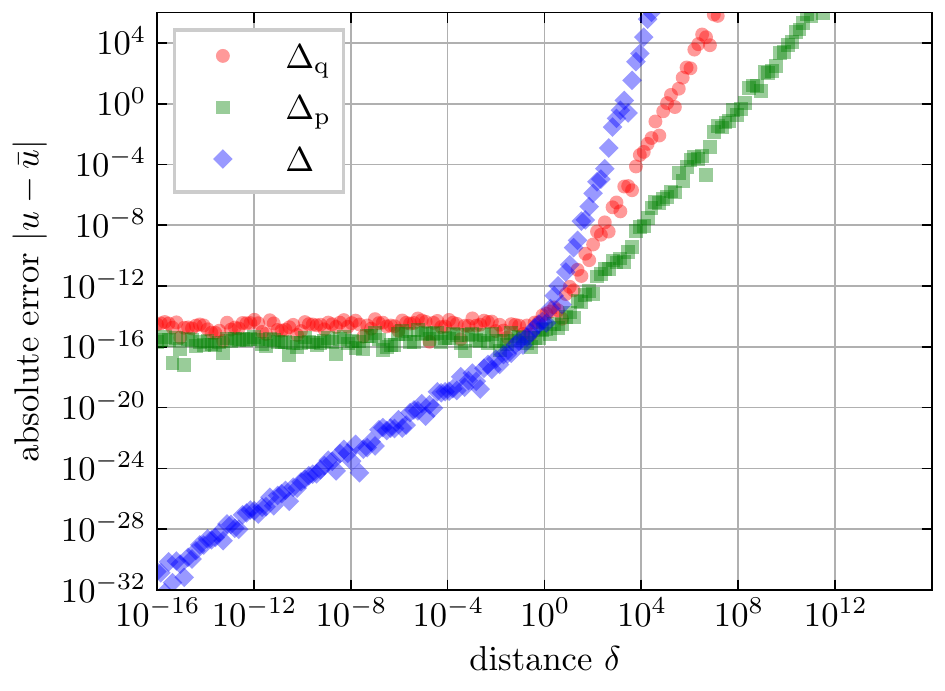}
    \caption{Critical case $\Delta \longrightarrow 0$, sum-of-products.}
    \label{fig:invariants-error-delta-sop}
    \end{subfigure}
    \hfill
    \begin{subfigure}[b]{0.45\textwidth}
    \centering
    \includegraphics[width=\textwidth]{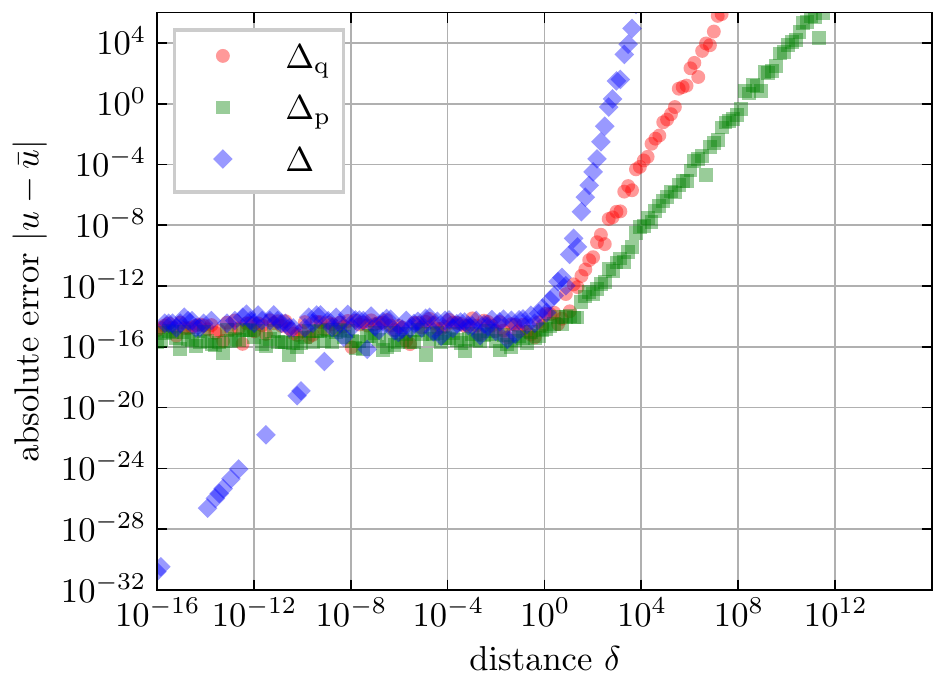}
    \caption{Critical case $\Delta \longrightarrow 0$, naive.}
    \label{fig:invariants-error-delta-naive}
    \end{subfigure}
    \vspace{1cm}

    \begin{subfigure}[b]{0.45\textwidth}
    \centering
    \includegraphics[width=\textwidth]{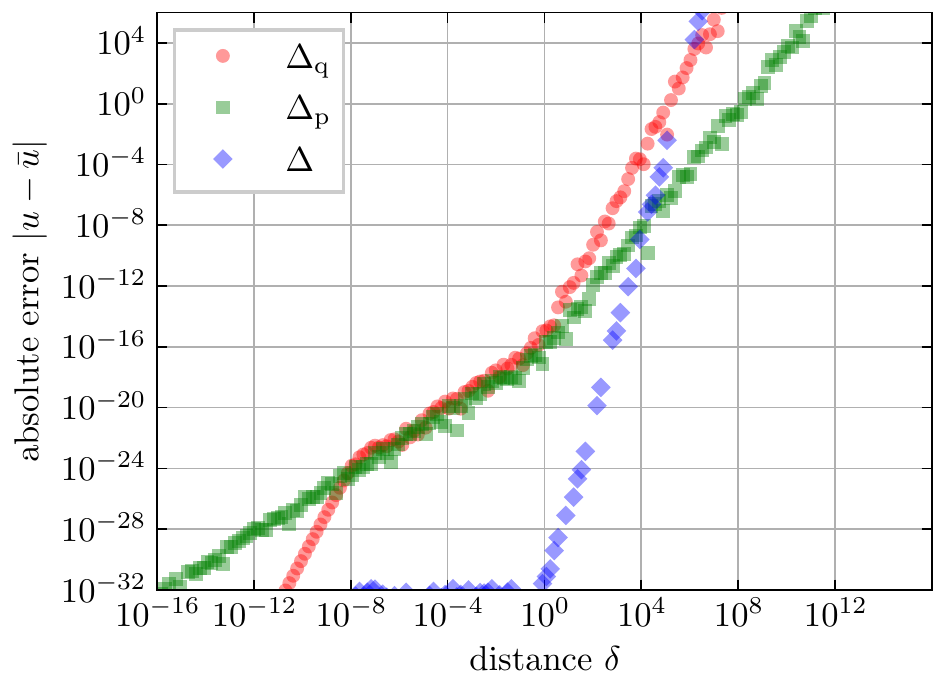}
    \caption{Critical case $\Delta_p \longrightarrow 0$, sum-of-products.}
    \label{fig:invariants-error-deltap-sop}
    \end{subfigure}
    \hfill
    \begin{subfigure}[b]{0.45\textwidth}
    \centering
    \includegraphics[width=\textwidth]{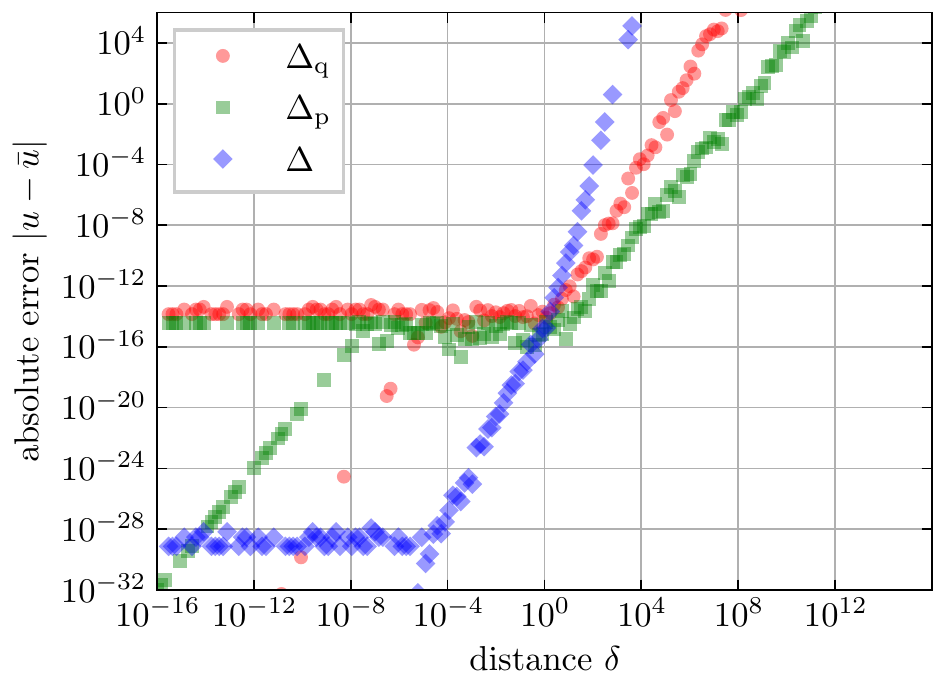}
    \caption{Critical case $\Delta_p \longrightarrow 0$, naive.}
    \label{fig:invariants-error-deltap-naive}
    \end{subfigure}
    \vspace{1cm}

    \begin{subfigure}[b]{0.45\textwidth}
    \centering
    \includegraphics[width=\textwidth]{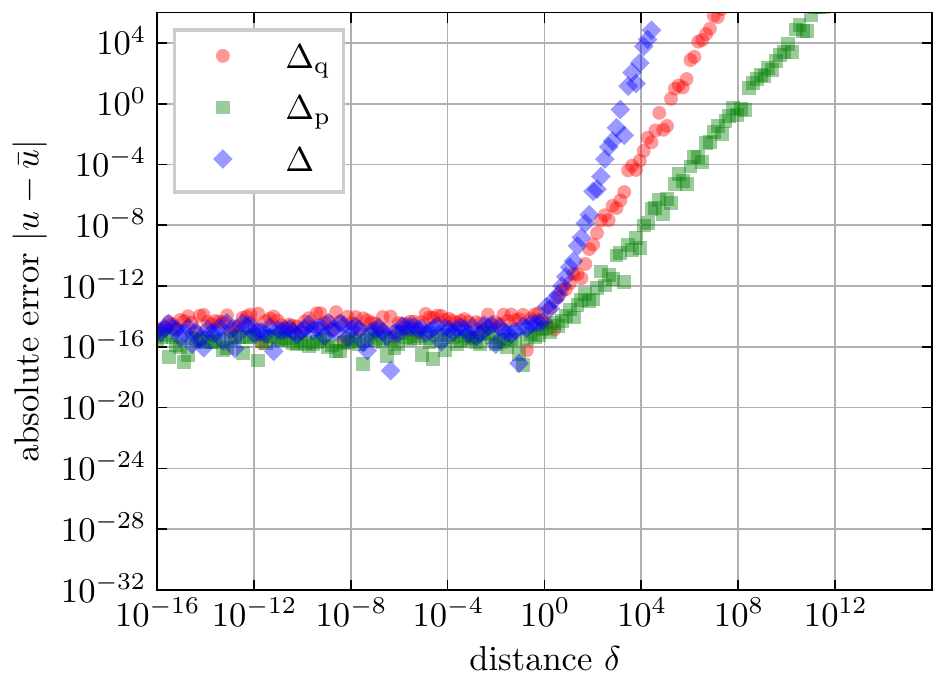}
    \caption{Critical case $\Delta_q \longrightarrow 0$, sum-of-products.}
    \label{fig:invariants-error-deltaq-sop}
    \end{subfigure}
    \hfill
    \begin{subfigure}[b]{0.45\textwidth}
    \centering
    \includegraphics[width=\textwidth]{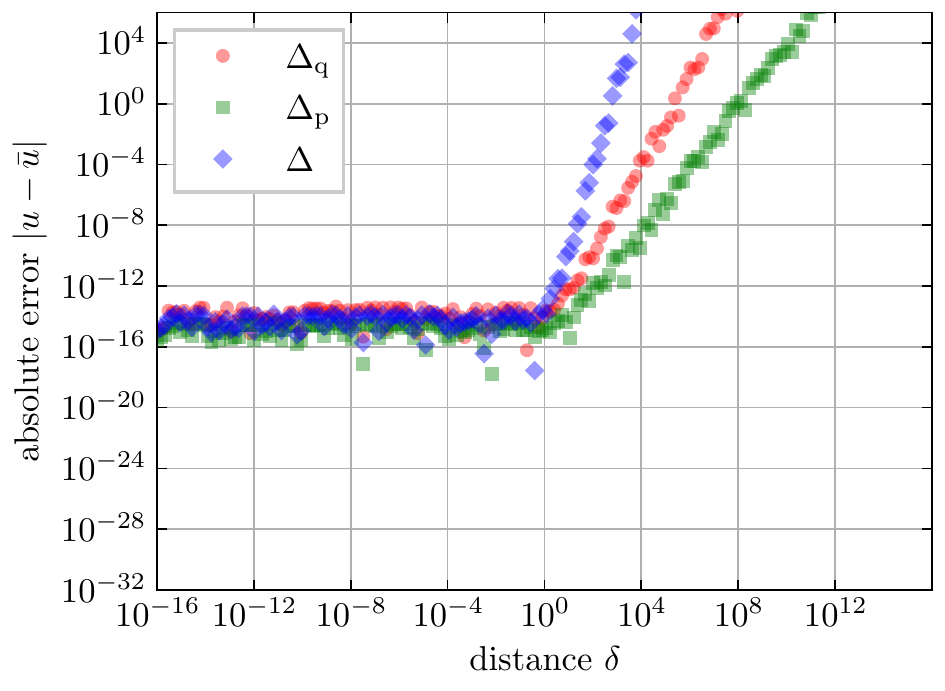}
    \caption{Critical case $\Delta_q \longrightarrow 0$, naive.}
    \label{fig:invariants-error-deltaq-naive}
    \end{subfigure}
    \vspace{1cm}

    \caption{Absolute error in matrix invariants for transformation matrix Case I.}
    \label{fig:invariants-error}
\end{figure}

Absolute errors in the computation of matrix invariants $\Delta_q, \Delta_p$ and $\Delta$ are 
shown in \autoref{fig:invariants-error}. 
Plots are produced for the transformation matrix Case I and various critical cases. 
Invariants computed using the sum-of-products formulas 
\eqref{eq:sop-discriminant-condensed3}, \eqref{eq:sop-deltap} and \eqref{eq:sop-deltaq} 
are included on the left while naive expressions \eqref{eq:naive-discriminant}, 
\eqref{eq:naive-deltap} and \eqref{eq:naive-deltaq} are used for the results plotted on the right.

The first critical case $\Delta \longrightarrow 0$ 
(\autoref{fig:invariants-error-delta-sop} and \autoref{fig:invariants-error-delta-naive}) 
shows decreased error for the discriminant (blue points) for distances $\delta < 1$. 
The naive expressions exhibit large rounding error and there is no measurable dependence  
for smaller $\delta$ distances. 
Note, that the absolute error should be decreasing, since the actual value of the discriminant 
$\Delta$ goes to zero. 
The sum-of-products expressions, however, consistently remain sensitive to distances $\delta$ 
as small as machine precision.

For the scenario of triple eigenvalues, as represented by the case $\Delta_p \longrightarrow 0$, 
and shown in 
\autoref{fig:invariants-error-deltap-sop} and \autoref{fig:invariants-error-deltap-naive}, the 
important difference is in the sensitivity of $\Delta_p$ (green squares) with respect to the 
distance $\delta$. 
Again, the sum-of-products approach shows improved accuracy where the absolute error in 
$\Delta_p$ decreases for smaller $\delta$. 

There is no significant difference between approaches for the last critical case 
$\Delta_q \longrightarrow 0$ 
(see \autoref{fig:invariants-error-deltaq-sop} and \autoref{fig:invariants-error-deltaq-naive})
since the (almost) sum-of-products expression for the $\Delta_q$ invariant is expected to suffer from 
subtractive cancellation between the terms $3 \Delta_{(0,1)(0,2)}$ and $4 \tr{\mathbf A} \Delta_p$.

\subsection{Eigenvalues}

\begin{figure}
    \centering
    \begin{subfigure}[b]{0.45\textwidth}
    \centering
    \includegraphics[width=\textwidth]{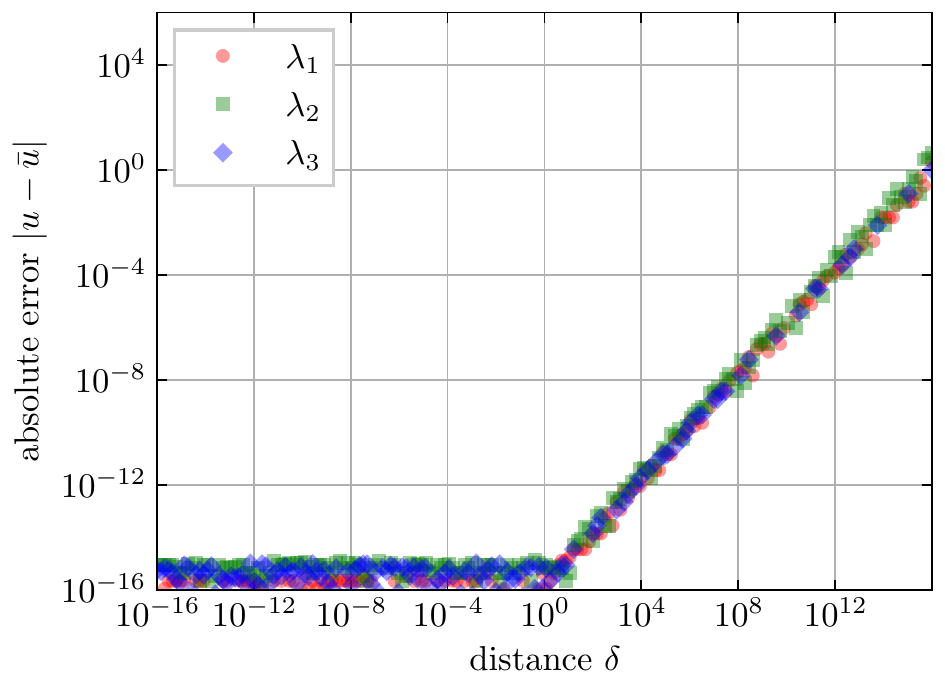}
    \caption{Critical case $\Delta \longrightarrow 0$, sum-of-products.}
    \label{fig:eigenvals-error-delta-sop}
    \end{subfigure}
    \hfill
    \begin{subfigure}[b]{0.45\textwidth}
    \centering
    \includegraphics[width=\textwidth]{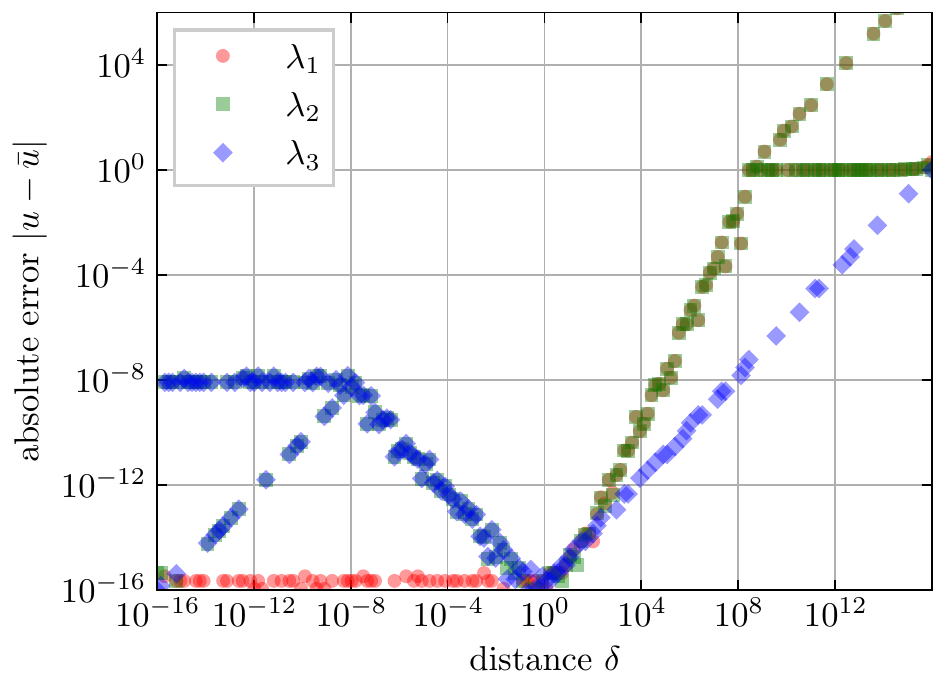}
    \caption{Critical case $\Delta \longrightarrow 0$, naive.}
    \label{fig:eigenvals-error-delta-naive}
    \end{subfigure}
    \vspace{1cm}

    \begin{subfigure}[b]{0.45\textwidth}
    \centering
    \includegraphics[width=\textwidth]{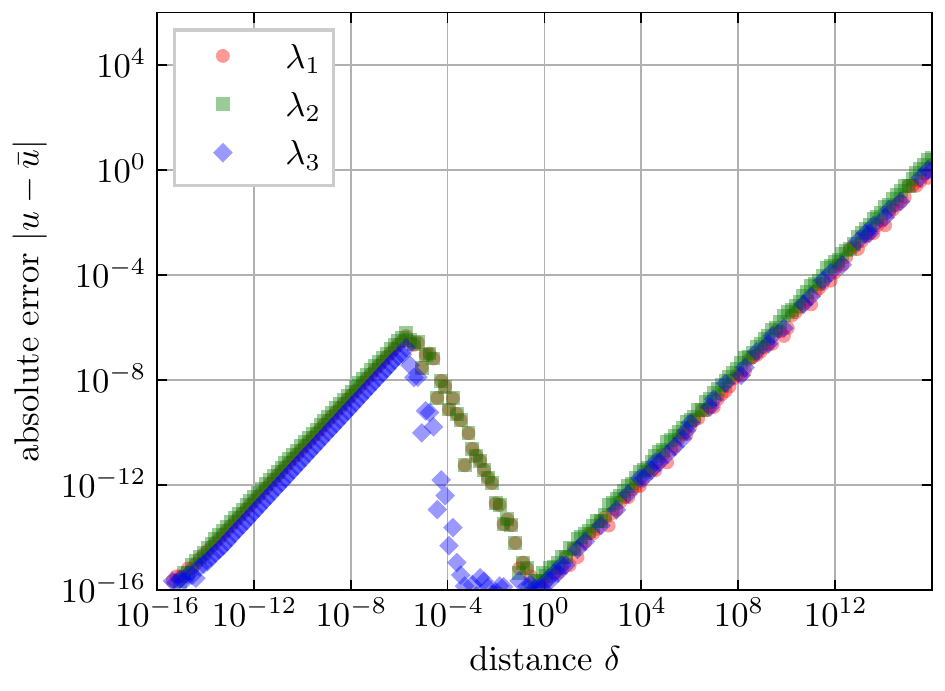}
    \caption{Critical case $\Delta_p \longrightarrow 0$, sum-of-products.}
    \label{fig:eigenvals-error-deltap-sop}
    \end{subfigure}
    \hfill
    \begin{subfigure}[b]{0.45\textwidth}
    \centering
    \includegraphics[width=\textwidth]{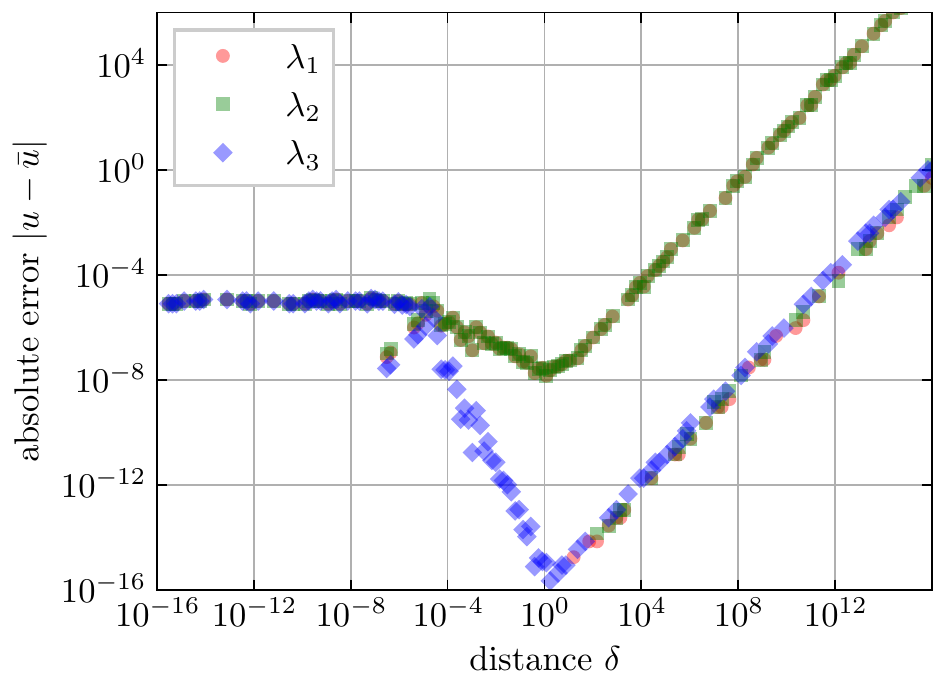}
    \caption{Critical case $\Delta_p \longrightarrow 0$, naive.}
    \label{fig:eigenvals-error-deltap-naive}
    \end{subfigure}
    \vspace{1cm}

    \begin{subfigure}[b]{0.45\textwidth}
    \centering
    \includegraphics[width=\textwidth]{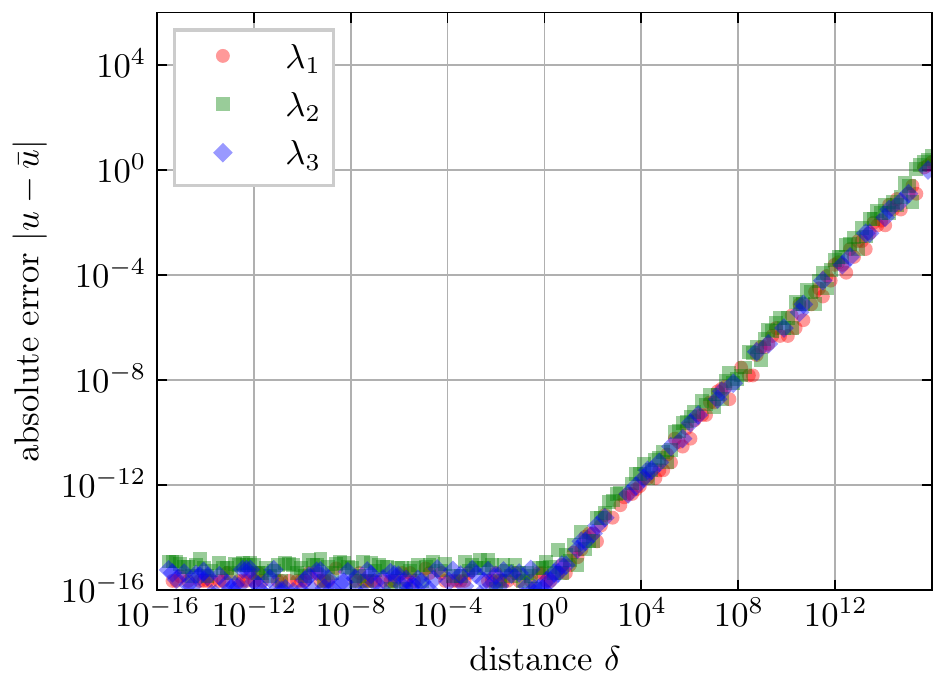}
    \caption{Critical case $\Delta_q \longrightarrow 0$, sum-of-products.}
    \label{fig:eigenvals-error-deltaq-sop}
    \end{subfigure}
    \hfill
    \begin{subfigure}[b]{0.45\textwidth}
    \centering
    \includegraphics[width=\textwidth]{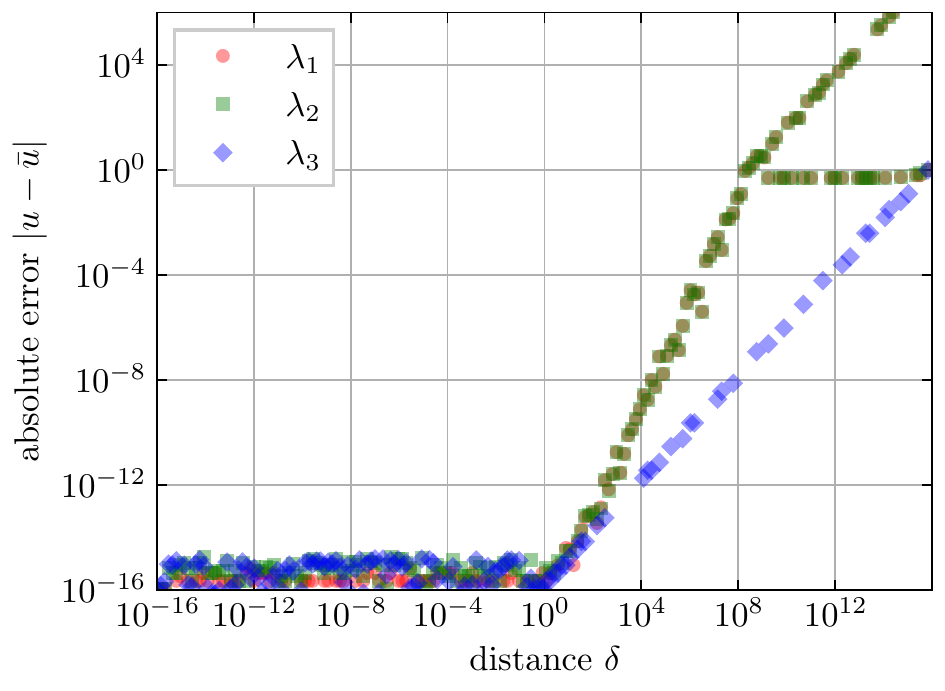}
    \caption{Critical case $\Delta_q \longrightarrow 0$, naive.}
    \label{fig:eigenvals-error-deltaq-naive}
    \end{subfigure}
    \vspace{1cm}

    \caption{Absolute error in eigenvalues for transformation matrix Case I.}
    \label{fig:eigenvalues-error}
\end{figure}

With formulas \eqref{eq:eigenvalues-dpdq} the absolute eigenvalues errors for the transformation matrix 
Case I are computed, see \autoref{fig:eigenvalues-error}. 
The error in each case is a direct consequence of the error shown for matrix invariants 
$\Delta, \Delta_p, \Delta_q$ and $I_1$. 
In general, accuracy is improved for sum-of-products expressions (left column).

For the critical case $\Delta_p \longrightarrow 0$ (\autoref{fig:eigenvals-error-deltap-sop} and 
\autoref{fig:eigenvals-error-deltap-naive}) there is an increase in the error around 
$\delta \approx \SI{e-6}{}$.
The reason for this error is contained in formula for the $\arctan$ argument in \eqref{eq:angleatan}, 
namely $\frac{3 \sqrt{3} \sqrt{\Delta}}{\Delta_q}$.
The discriminant $\Delta$ should be equal to zero (in infinite precision arithmetic, two smallest 
eigenvalues are equal), but rounding causes the discriminant to be nonzero 
(see e.g. \autoref{fig:invariants-error-deltap-sop}). 
The floating point value of the discriminant is dominated by rounding error, so it has no significance 
(the relative error is infinitely large). 
As $\Delta_q$ becomes sufficiently small it amplifies insignificant values in the discriminant and 
increases the error of the fraction.

\paragraph{Remark:} %
Absolute errors in \autoref{fig:invariants-error} and \autoref{fig:eigenvalues-error} were computed for 
the transformation matrix $\mathbf U$, Case I, see \autoref{tab:similarity-trans-bench}. 
Case II contains the parameter $\gamma > 0$ which controls the conditioning number $\kappa(\mathbf U)$ 
of the transformation matrix (the smaller $\gamma$, the larger the conditioning number). 
For the critical case $\Delta \longrightarrow 0$ the test matrix $\mathbf B$ is expressed as a function 
of the distance $\delta$ and the parameter $\gamma$,
\begin{align}
    \mathbf B = \mathbf U(\gamma) \mathbf \Lambda(\delta) \mathbf{U(\gamma)}^{-1} =
    \begin{bmatrix}
        \frac{- \delta + \gamma - 2}{\gamma} & - \frac{\delta + 2 \gamma + 2}{\gamma} & \frac{\delta + 2}{\gamma}\\- \frac{\delta + 2}{\gamma} & - \frac{\delta + \gamma + 2}{\gamma} & \frac{\delta + 2}{\gamma}\\\frac{\gamma - \left(\delta + 1\right) \left(\gamma + 2\right) - 2}{\gamma} & - \delta - \frac{2 \delta}{\gamma} - 4 - \frac{4}{\gamma} & \frac{\left(\delta + 1\right) \left(\gamma + 2\right) + 2}{\gamma}
    \end{bmatrix}.
\end{align}
The important observation is that for $\gamma \longrightarrow 0$ every component in the matrix goes to 
$\pm \infty$ with increasing proportionality coefficients in $\delta$. 
In formulas \eqref{eq:sop-discriminant-condensed3} the components of the vectors $\mathbf{\bar x}$ and 
$\mathbf{\bar y}$ could happen to be close to zero due to cancellation in products of matrix components. 
With a badly conditioned transformation matrix $\mathbf U$ the error in the cancellation becomes larger. 
This scaling affects the naive formula for discriminant \eqref{eq:naive-discriminant} negatively too.

\begin{figure}
    \centering
    \begin{subfigure}[b]{0.45\textwidth}
    \centering
    \includegraphics[width=\textwidth]{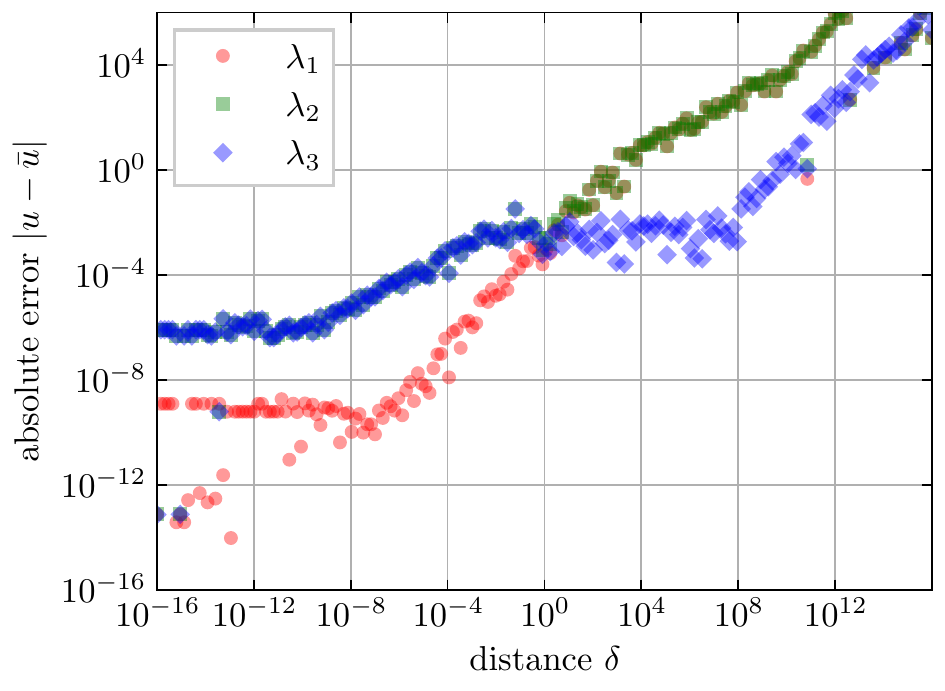}
    \caption{Transformation matrix Case II, sum-of-products.}
    \end{subfigure}
    \hfill
    \begin{subfigure}[b]{0.45\textwidth}
    \centering
    \includegraphics[width=\textwidth]{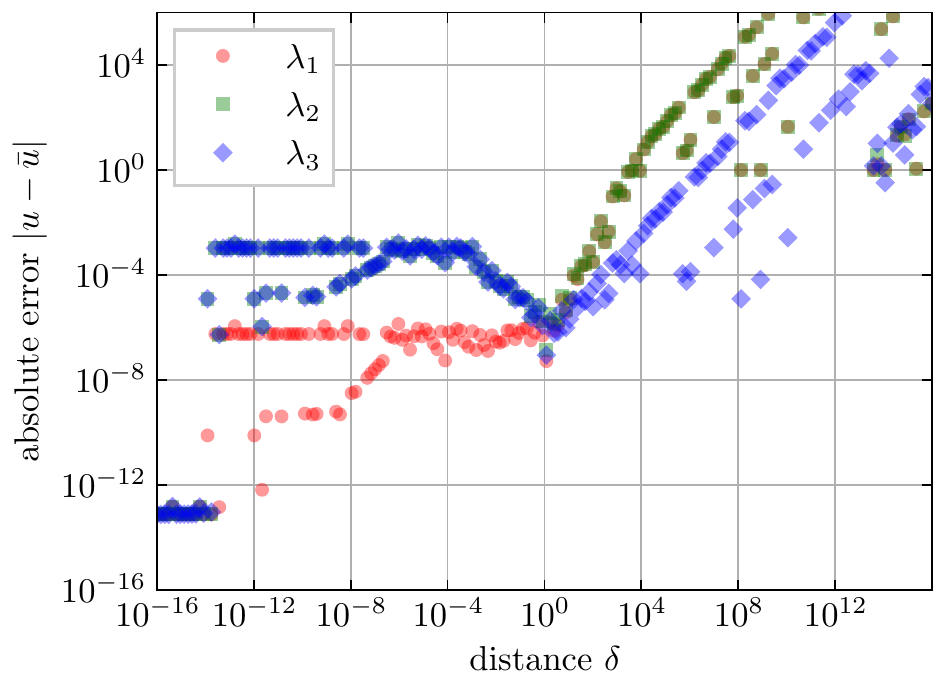}
    \caption{Transformation matrix Case II, naive.}
    \end{subfigure}
    \vspace{1cm}

    \caption{Absolute error in eigenvalues for critical case $\Delta \longrightarrow 0$ and 
             badly conditioned transformation matrix $\mathbf U$ with $\gamma = \SI{e-3}{}$.}
    \label{fig:eigenvalues-error-caseii}
\end{figure}

The effect of a badly conditioned matrix $\mathbf U$ with $\gamma = \SI{e-3}{}$ is depicted in 
\autoref{fig:eigenvalues-error-caseii}. 
Both, sum-of-products and naive expressions, show decreased accuracy for all eigenvalues.

\subsection{Eigenprojectors}

\begin{figure}
    \centering
    \begin{subfigure}[b]{0.45\textwidth}
    \centering
    \includegraphics[width=\textwidth]{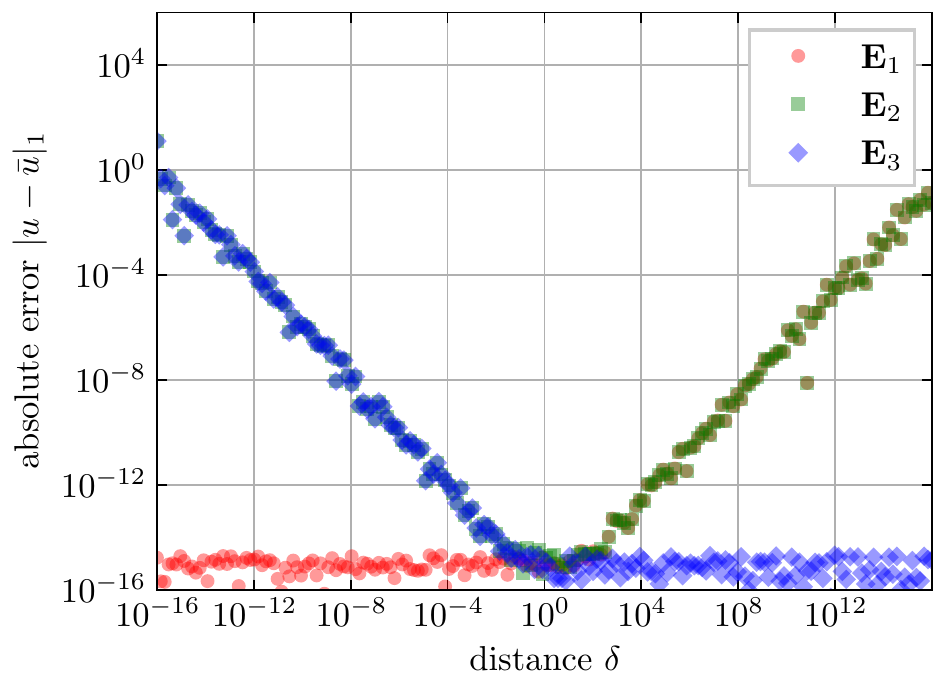}
    \caption{Critical case $\Delta \longrightarrow 0$, sum-of-products.}
    \end{subfigure}
    \hfill
    \begin{subfigure}[b]{0.45\textwidth}
    \centering
    \includegraphics[width=\textwidth]{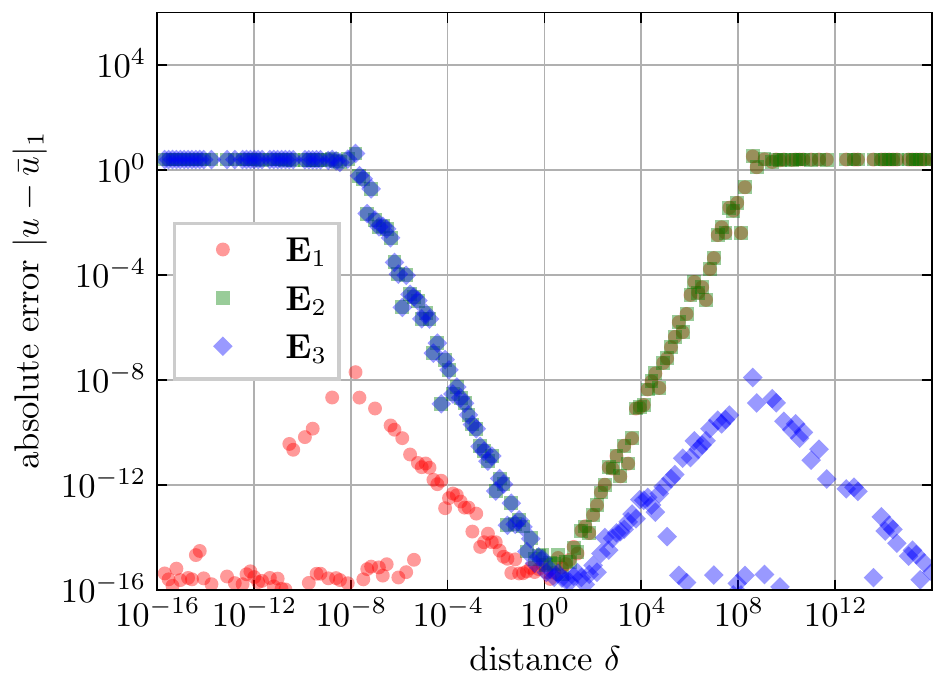}
    \caption{Critical case $\Delta \longrightarrow 0$, naive.}
    \end{subfigure}
    \vspace{1cm}

    \caption{Absolute $l_1$ error in eigenprojectors for transformation matrix Case I.}
    \label{fig:eigenproj-error}
\end{figure}

The eigenprojectors can be computed as eigenvalues derivatives, see \eqref{eq:eigenprojector}. 
The observed absolute error of this approach is outlined in \autoref{fig:eigenproj-error}. 
This error shows different characteristics as the absolute error observed for invariants and eigenvalues. 
The computed eigenprojectors satisfy the "normalisation properties" according to 
\eqref{eq:eigenprojectorsproperties}, which makes their absolute error resemble more the relative error 
of eigenvalues.

\subsection{Matrix powers}

A typical use case for symbolic spectral decomposition is the evaluation of matrix functions $f(\mathbf A)$ 
according to Sylvester interpolating definition (see \citep{higham2008functions})
\begin{align}
    f(\mathbf A) = \sum\limits_{k=1}^3 f(\lambda_k) \mathbf E_k.
\end{align}

\begin{figure}
    \centering
    \begin{subfigure}[b]{0.45\textwidth}
    \centering
    \includegraphics[width=\textwidth]{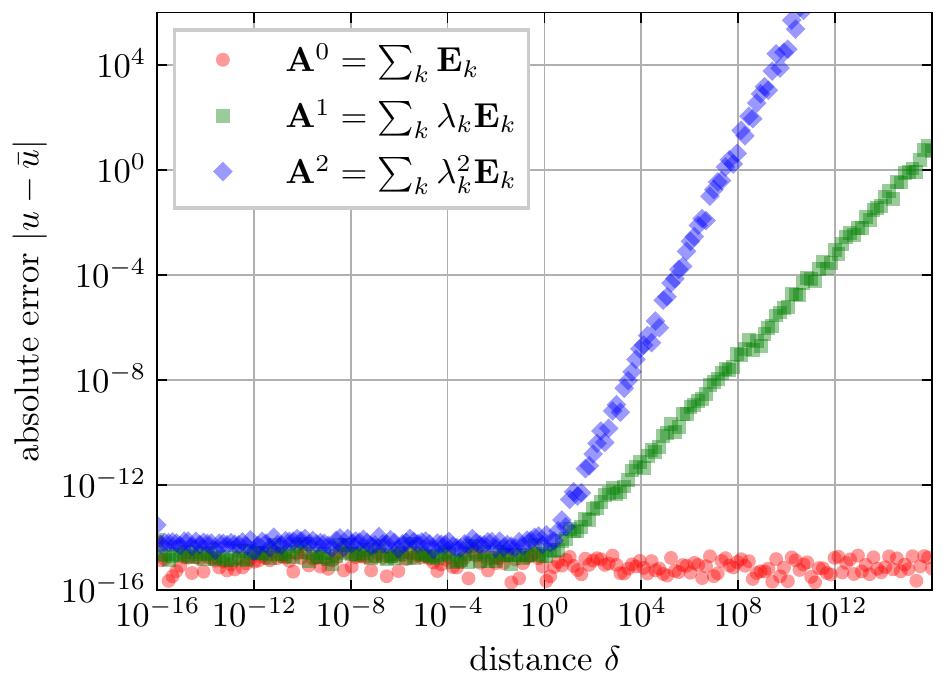}
    \caption{Critical case $\Delta \longrightarrow 0$, sum-of-products.}
    \label{fig:eigpow-sop-error}
    \end{subfigure}
    \hfill
    \begin{subfigure}[b]{0.45\textwidth}
    \centering
    \includegraphics[width=\textwidth]{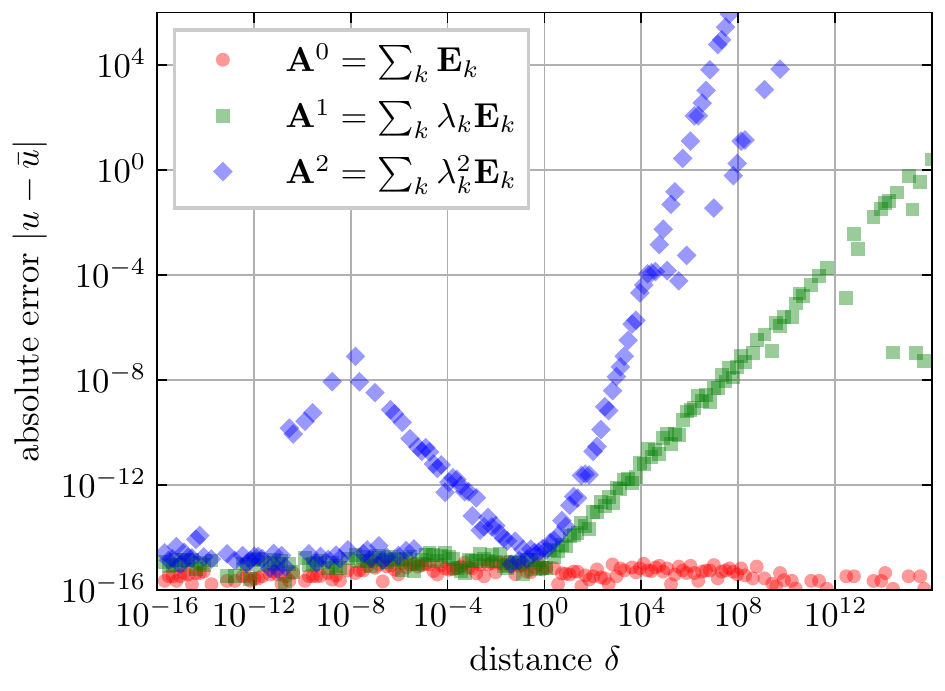}
    \caption{Critical case $\Delta \longrightarrow 0$, naive.}
    \label{fig:eigpow-naive-error}
    \end{subfigure}
    \vspace{1cm}

    \caption{Absolute $l_1$ error in matrix powers.}
    \label{fig:eigenpow-error}
\end{figure}

The upper bound for the error in this formula is then a consequence of errors in eigenvalues 
and eigenprojectors (due to trivial triangle inequality estimate). 
Benchmark results for the critical case $\Delta \longrightarrow 0$ and a well conditioned 
transformation matrix (Case I) are included in \autoref{fig:eigenpow-error}. 
In \autoref{fig:eigpow-sop-error} and \autoref{fig:eigpow-naive-error} the main contribution to the 
error is due to eigenvalues errors, which is then scaled with a factor corresponding to the matrix power. 
There is a notable improvement in the accuracy for $\mathbf A^2$ when the sum-of-products approach is used.

\section{Conclusion}

This paper studies the closed-form solution approach to spectral decomposition of real-valued 
matrices with real eigenvalues. 
Formulas based on trigonometric transformations are derived in the beginning. 
The notion of the matrix discriminant is then generalised based on the determinant of certain matrix powers. 
The special case of $3 \times 3$ matrices is studied for associated rounding errors in the symbolic computation 
of the spectral decomposition and several potential sources of catastrophic cancellation are identified. 
Discriminant and sub-discriminants are expressed on the basis of sum-of-products expressions. 
It is shown, that this mathematically equivalent procedure provides alternative expressions for particular 
matrix invariants (discriminant and invariant $\Delta_p$), which do not suffer catastrophic cancellation and 
consequently show much improved numerical floating point accuracy in all critical cases.

A set of numerical benchmarks is executed to test newly developed expressions. 
Absolute errors in matrix invariants $\Delta, \Delta_p$ and $\Delta_q$ show different characteristics and 
are in general superior to the naive expressions used in existing literature. 
With the employed techniques the error in computed eigenvalues can be decreased from half to full 
machine precision.

The overall improved algorithm is included in \autoref{sec:algorithm} and is believed to serve as basis 
for those who require near machine precision of eigenvalues and preserved symbolic, functional dependence 
of the eigendecomposition on the matrix elements.

An extension of the presented formalism to matrices having a spectral decomposition over the complex numbers
or general complex-valued matrices is considered straightforward.
While the advantages of the sum-of-products based (sub)-discriminant(s) persist, the main difference lies then
in an appropriate choice for Equation \eqref{eq:charpolysubst} with trivial subsequent steps.

\clearpage

\appendix
\section{Improved algorithm for spectral decomposition}
\label{sec:algorithm}

Following algorithm is written in Python programming language and tested using \textbf{Python 3.9.7} and symbolic algebra package \textbf{SymPy 1.8}.

\begin{small}
\begin{minted}{python}
import sympy as sy
import itertools


def subdiscriminant(A, subU, subV):
    U = sy.matrices.Matrix.vstack(*[A.pow(p).vec().T for p in subU])
    V = sy.matrices.Matrix.hstack(*[A.pow(p).T.vec() for p in subV])

    enums = list(itertools.combinations(range(U.cols), U.rows))
    u = [U[:, e].det() for e in enums]
    v = [V[e, :].det() for e in enums]
    return sum([uu * vv for uu, vv in zip(u, v)])


def eigenvalues(A):
    I1 = sy.Trace(A)
    dp = 1 / 2 * subdiscriminant(A, [0, 1], [0, 1])
    dq = 3 * subdiscriminant(A, [0, 1], [0, 2]) - 4 * I1 * dp
    disc = subdiscriminant(A, [0, 1, 2], [0, 1, 2])
    atan_arg = 3 * sy.sqrt(3) * sy.sqrt(disc) / dq
    phi = sy.atan(atan_arg) + (1 - sy.sign(dq)) * sy.pi / 2

    return [(I1 + 2 * sy.sqrt(dp) * sy.cos((phi + 2 * sy.pi * k) / 3)) / 3
            for k in [1, 2, 3]]


def eigenprojectors(A):
    eigenvals = eigenvalues(A)
    return [eigenval.diff(A).T for eigenval in eigenvals]
\end{minted}
\end{small}

\bibliographystyle{abbrvnat}
\bibliography{references}


\end{document}